\documentclass[sn-mathphys,Numbered]{sn-jnl}

\usepackage{multirow}%
\usepackage{amsmath,amssymb,amsfonts}%
\usepackage{amsthm}%
\usepackage{mathrsfs}%
\usepackage[title]{appendix}%
\usepackage{xcolor}%
\usepackage{textcomp}%
\usepackage{manyfoot}%
\usepackage{booktabs}%
\usepackage{algorithm}%
\usepackage{algorithmicx}%
\usepackage{algpseudocode}%
\usepackage{listings}%
\usepackage{moreverb,lineno}
\usepackage{enumerate}
\usepackage{overpic}
\usepackage{graphicx}
\newcommand\BibTeX{{\rmfamily B\kern-.05em \textsc{i\kern-.025em b}\kern-.08em
		T\kern-.1667em\lower.7ex\hbox{E}\kern-.125emX}}

\usepackage{setspace}
\if@secthm
\newtheorem{theorem}{Theorem}[section]
\else

\fi
\usepackage{amsmath,amsthm,amssymb,dsfont}
\def\cals{\mathcal{S}}
\def\call{\mathcal{L}}
\def\calu{\mathcal{U}}
\def\calk{\mathcal{K}}

\def\calo{\mathcal{O}}
\def\M{\cals}
\def\R{\mathbb{R}}
\def\N{\mathbb{N}}
\def\C{\mathbb{C}}
\def\bigdot{\boldsymbol{\cdot}}
\def\eref#1{{\rm (\ref{#1})}}
\def\n{\mathbf{n}}

\def\to{\rightarrow}
\def\cp{{\mathrm{cp}}}
\DeclareMathOperator*{\arginf}{arg\,inf}
\DeclareMathOperator*{\esssup}{ess\,sup}

\def\Span{\mbox{span}}

\usepackage{lmodern}
\usepackage[T1]{fontenc}
\renewcommand\rmdefault{ptm} 
\DeclareSymbolFont{myletters}{OML}{ztmcm}{m}{it}
\DeclareMathSymbol{\uplambda}{\mathord}{myletters}{"15}

\def\dS{{d_\M}}
\def\TM{\mathcal{T}_x{\M}}

\def\dS{{d_\M}}

\usepackage{multirow}

\usepackage{color}

\usepackage{titlesec}
\newtheorem{example}[subsection]{Example}
\def\exmp#1#2{{\begin{example}\label{#2}\textbf{#1.}\end{example}\noindent}}

\usepackage[normalem]{ulem}
\usepackage{color} 

\newcommand{\blue}[1]{{\color{blue}{#1}}}

\definecolor{mygray}{gray}{0.4}
\definecolor{darkgreen}{rgb}{0, 0.45, 0}

\newcommand{\green}[1]{{\color{black} #1}}

\newcommand{\teal}[1]{{\color{black} #1}}
\usepackage{soul}
\usepackage[normalem]{ulem}
\newcommand{\CMadd}[1]{\green{#1}}    
\newcommand{\CMout}[1]{\teal{\sout{#1}}}    
\newcommand{\CMrep}[2]{\teal{#2}} 

\begin{document}




\title{Exploring Oversampling in RBF Least-Squares Collocation Method of Lines for Surface Diffusion}


\author*[1,2]{\fnm{Meng} \sur{Chen}}\email{chenmeng{\_}math@ncu.edu.cn}

\author[3]{\fnm{Leevan} \sur{Ling}}\email{lling@hkbu.edu.hk}

\affil*[1]{\orgdiv{Department of Mathematics}, \orgname{Nanchang University}, \orgaddress{\street{Xuefu Avenue}, \city{Nanchang}, \postcode{330031}, \state{Jiangxi}, \country{China}}}

\affil[2]{\orgdiv{Institute of Mathematics and Interdisciplinary Sciences}, \orgname{Nanchang University}, \orgaddress{\street{Xuefu Avenue}, \city{Nanchang}, \postcode{330031}, \state{Jiangxi}, \country{China}}}

\affil[3]{\orgdiv{Department of Mathematics}, \orgname{Hong Kong Baptist University}, \orgaddress{\street{Waterloo Road}, \state{Kowloon Tong}, \country{Hong Kong}}}

\abstract{This paper investigates the numerical behavior of the radial basis functions least-squares collocation (RBF-LSC) method of lines (MoL) for solving surface diffusion problems, building upon the theoretical analysis presented in [\href{https://doi.org/10.1137/21M1444369}{SIAM J. Numer. Anal., 61 (3), 1386-1404}]. Specifically, we examine the impact of the oversampling ratio, defined as the number of collocation points used over the number of RBF centers for quasi-uniform sets, on the stability of the eigenvalues, time stepping sizes taken by Runge-Kutta methods, and overall accuracy of the method. By providing numerical evidence and insights, we demonstrate the importance of the oversampling ratio for achieving accurate and efficient solutions with the RBF-LSC-MoL method. Our results reveal that the oversampling ratio plays a critical role in determining the stability of the eigenvalues, and we provide guidelines for selecting an optimal oversampling ratio that balances accuracy and computational efficiency.
}

\keywords{Surface diffusion, kernel-based least-squares collocation methods, method of lines, partial differential equations on manifolds, eigenvalue stability.}



\maketitle
\newpage
\section{Introduction}
\green{
 The Radial Basis Function (RBF) methods have emerged as promising numerical method for solving partial differential equations (PDEs) due to their mesh-free nature and inherent capability to handle complex geometries, such as surfaces. The concept of overtesting is to use more collocation conditions than the number of basis functions. This results in an overdetermined matrix system after spatial discretization. Specifically, the Least-Squares Collocation (LSC) method is derived when the minimum $\ell^2$-norm solution is selected.
 Considerable progress has been made in the development of various LSC methods for elliptic PDEs in bounded domains \cite{Cheung+LingETAL-leaskerncollmeth:18}, and on surfaces \cite{Chen+Ling-Extrmeshcollmeth:20,Cheung+Ling-Kernembemethconv:18}. These methods have been proven to be convergent, and interestingly, theoretical error bounds hold under conditions of a sufficiently high oversampling ratio. According to numerical evidence, a high oversampling ratio does not negatively impact accuracy, but computational complexity is linearly related to the oversampling ratio.

 In the context of parabolic PDEs, the Method of Lines (MoL) firstly discretizes in space, suggesting that RBF methods may also benefit from oversampling. Our recent study \cite{Chen+CheungETAL-kernleascollmeth:23} showed that overtesting is also a theoretical requirement for achieving the error bound of an RBF-LSC-MoL numerical solution for parabolic surface diffusion PDEs.
Technically, adjusting the oversampling ratio means employing a different spatial discretization, and it influences the properties of the ODEs in MoL. This distinction from the elliptic case, where the effects are seen in accuracy and convergence rates, is the central focus of this paper.

This paper is structured as follows: Section~\ref{sec:Notations and preliminaries} begins by reviewing the necessary assumptions and definitions relevant to surface diffusion problems. Section~\ref{SecRBFMoL} presents the RBF-LSC-MoL analyzed in \cite{Chen+CheungETAL-kernleascollmeth:23}. The rest of the paper, detailed in Section~\ref{sec:Eigenvalue stability}, employs a series of examples to explore the role of the oversampling ratio in influencing the stability of the eigenvalues. We start with a linear isotropic diffusion-reaction problem, progress to an anisotropic case with a space-dependent diffusion tensor, and finally simulate two surface diffusion problems without fine-tuning of setups. All these numerical evidence are to provide readers some insights in picking the oversampling ratio for RBF-LSC-MoL.
}

\section{Surface diffusion: notations and preliminaries}\label{sec:Notations and preliminaries}
We consider a closed, connected, orientable, and complete Riemannian manifold  $\M\subset\R^d$
of dimension $\text{dim}\,\M := \dS= d-1$ and  of class $\mathcal{C}^{m+1}$, $m\in\N$, where the theories in \cite{Maerz+Macdonald-CALCSURFWITHGENE:12,Hangelbroek+NarcowichETAL-DireInveResuBoun:17} are applicable.
Our interest lies in solving parabolic PDEs for some time-dependent surface scalar function $u:\M\times[0,T]\to\R$  in the form of
\begin{subequations}\label{eqPDE}
	\begin{align}
		\displaystyle
		\dot u(x,t) + \call_\M u(x,t) = f(x,t) &\quad \mbox{for }(x,t)\in\M\times[0,T], \label{eqPDE-a}
		\\
		u(x,0) = g(x)~~~\, &\quad \mbox{for } x\in\M. \label{eqPDE-b}
	\end{align}
\end{subequations}
for some given square integrable functions $f:\M\times[0,t]\to \R$ and $g:\M\to \R$, and  a second-order uniformly elliptic operator   in divergence form
\begin{equation}\label{equLs}
	\call_{\M } u(x,t):= - \nabla_{\M }\bigdot \big( A(x,t) \nabla_{\M } u(x,t) \big) +b\,u(x,t) \quad \mbox{for }(x,t)\in\M\times[0,T],
\end{equation}
where  $b$ is a constant, and  the diffusion tensor $A$ intrinsic to surface$:\M\times[0,t]\to \R^{d\times d}$  satisfies,
for all $x\in\M$ and $t\in[0,T]$, some symmetric positive definiteness assumption \cite[Asm. 1]{Chen+CheungETAL-kernleascollmeth:23} when restricted to the tangent space $\TM$ of $\M$ at $x\in\M$.
The surface diffusion operator in \eref{equLs} is defined by a smooth vector field $\n :\M\to \R^d$  that spans the normal space and the $\TM$-orthogonal projection matrix
\begin{equation}\label{def:P(y)}
	P(x) := I_d - \n(x) \n(x)^T.
\end{equation}
Then, the surface gradient operator $\nabla_\M$
\begin{equation}\label{def:grad}
	\nabla_\M v(x) := P(x)\nabla\big(v\circ\cp\big)(x)  = \nabla\big(v\circ\cp\big)(x),\quad x\in\M,
\end{equation}
and surface divergence operator $\nabla_\M\bigdot$
\begin{equation}\label{def:div}
	\nabla_\M \bigdot \mathbf{g}(x) := \big(P(x) \nabla\big)\bigdot\big(\mathbf{g}\circ\cp\big)(x) =  \nabla\bigdot\big(\mathbf{g}\circ\cp\big)(x),\quad x\in\M,
\end{equation}
are defined, respectively, for  continuously differentiable scalar surface function $v:\M\to\R$ and for continuously differentiable surface  vector field $\mathbf{g}:\M\to\R^{d}$. In \eref{def:grad} and \eref{def:div}, $\cp$ is the $\mathcal{C}^{m}$-smooth Euclidean closest point retraction map:
\begin{equation}\label{def:cp}
	\cp(y) := \arginf_{x\in\M} \|x-y\|_{L^2(\R^d)}:\Omega\to \M,
\end{equation}
in some narrow-band domain $\Omega$ containing $\M$.
\green{The introduction of $\cp$ enables functions defined on the surface to be evaluated off the surface. As a result, the standard differential operators are well-defined in \eqref{def:grad} and \eqref{def:div}.}

\section{RBF-collocation method of lines}\label{SecRBFMoL}
In \cite{2005Approximation,2020Solving,2012Solving}, spherical RBF approximation methods with the method of lines (MoL) were applied to solving parabolic PDEs on unit spheres. Backward Euler and Crank-Nicolson schemes were respectively employed in \cite{2005Approximation} with meshless collocation approaches on fixed surfaces, and in \cite{2012Solving} with collocation and other two (Galerkin and RBF-FD) methods on evolving spheres. Meanwhile, the authors of \cite{2020Solving} used Laplace transforms and quadrature in time along with Galerkin approximation in space.

We will not go into detail about  the semi-discretized solution and instead present the fully discretized RBF least-squares collocation method of lines (RBF-LSC-MoL) \CMrep{in}{from} \cite{Chen+CheungETAL-kernleascollmeth:23} \CMrep{with the}{using a} minimal set of notations. First, we pick a  global, symmetric positive definite, and Sobolev space $H^{m+1/2}(\R^d)$ reproducing \cite{Matern} kernels
$\Psi (\cdot,\cdot)= {\Psi}_{m+1/2}(\cdot,\cdot): \R^d \times \R^d \rightarrow \R$\CMout{to $\M$},
whose Fourier transforms with $\|{\omega}\|_2=\text{dist}(\cdot,\cdot)$  decay like
\begin{equation}\label{kernelFour}
 	c_1(1+\|{\omega}\|_2^2)^{-{(m+1/2)}}\le  \hat{\Psi}({\omega})
 	\leq c_2(1+\|{\omega}\|_2^2)^{-{(m+1/2)}} \quad \mbox{for all }{\omega}\in\R^d,
 \end{equation}
for some constants $0<c_1\leq c_2$; for example, we can choose between two types of commonly-used kernels:  the standard Whittle-Mat\'{e}rn-Sobolev kernels \cite{Matern} or the  Wendland compactly supported kernels \cite{Wendlandfun} with a smoothness order $m+1/2$.
The extra half order of smoothness in $\Psi$ ensures that the corresponding restricted surface kernel  reproduces  $H^{m}(\M)$ \cite{Fuselier,narcowich2007approximation}.

For some set of $n_Z$ trial centers $Z=\{z_1, \ldots, z_{n_Z}\} \subset \M$ with
fill distance $h_{Z}$ and the {separation distance} $q_{Z}$
\begin{equation}\label{def_h_q}
 	h_{Z}:=\sup_{\zeta\in \M} \inf_{\eta\in Z} \text{dist}({\zeta},{\eta}) \quad \mbox{and}\quad q_{Z}:=\frac{1}{2}\inf_{{{z}_i\neq {z}_j \in Z}}\text{dist}({z}_i,{z}_j),
 \end{equation}
we define the RBF trial space by
\begin{equation}\label{Uzfinite}
 	\calu_{Z,\M ,\Psi}:= \Span\{\Psi\left(\cdot\, ,z_j\right)\;|\;z_j\in Z\}.
\end{equation}
We assume that $Z$ is quasi-uniform with  \teal{a} bounded {mesh ratio}  $\rho_{Z} =h_Z/q_Z \geq 1$ as $Z$ is \CMout{being} refined.
Our goal is to find a time-dependent numerical solution to the surface diffusion system \eref{eqPDE} in the trial space using a set of unknown coefficient functions of time
$\lambda_Z(t):=[\lambda_1(t),\ldots,\lambda_{n_Z}(t)]^T$, that is,
\begin{equation}\label{eq:trial sp0}
 	u_Z (\cdot, t)   = \sum_{z_j\in Z}^{n_Z}\lambda_j(t)  \Psi(\cdot,z_j)
 	=: \Psi( \cdot, Z) \lambda_Z(t) .
 \end{equation}
By putting the ansatz \eref{eq:trial sp0} into the governing equation \eref{eqPDE-a} and oversampling by collocation at some sufficiently dense set of quasi-uniform  points $X = \{x_1,\ldots, x_{n_X}\} \subset \M$ that  satisfies the following the denseness  constraints \cite[Lem. 3.5]{Chen+CheungETAL-kernleascollmeth:23}
\begin{equation}\label{eq:desnseness stability}
 	\left\{
 	\begin{array}{rl}
 		C h_X^{2m- 4}  h_Z^{-2m}  < \frac14 & \mbox{if }h_Z \leq 1, \\
 		C h_X^{2m- 4}  h_Z^{-2m+4} < \frac14& \mbox{otherwise},\\
 	\end{array}
 	\right.
 \end{equation}
for some $C$  depending on ${\Psi},  \M, \call, \rho_Z$ and  $\rho_X$, we obtain the following
\teal{semi-discretized} overdetermined differential  algebraic equations (DAE)
\begin{equation}\label{eq:DAE}
 	\big[\Psi( X, Z)\big] \dot \lambda_Z(t)  + \big[\call_\M  \Psi( X, Z)\big] \lambda_Z(t) = f(X,t),
 \end{equation}
where the $ij$-entry of the $n_X\times n_Z$ matrix $[\Psi( X, Z)]$  (and $[\call_\M\Psi( X, Z)]$) is defined by $\Psi(x_i,z_j)$  (and $\call_\M\Psi(x_i,z_j)$) for $x_i\in X$ and $z_j\in Z$.  Instead of coefficient functions, one may \teal{prefer} to work with nodal values $u_Z(t)$ of $u$ at $Z$ by considering
\begin{equation}\label{eq:DAE2}
 	\big[\Psi( X, Z) \Psi(Z,Z)^{-1}\big] \dot  u_Z(t)  +  \big[\call_\M  \Psi( X, Z)   \Psi(Z,Z)^{-1}\big] u_Z(t) = f(X,t).
\end{equation}
\green{
A semi-discretized solution $u_{Z,\alpha}$  can then be defined by the time-integral, i.e., the $L^2[0,T]$-norm, of the $\ell^2$-residual at $X$ of the DAE in \eqref{eq:DAE}, subject to appropriate initial conditions. }
We have the following estimate in \cite[Thm. 3.6 \& Cor. 4.2]{Chen+CheungETAL-kernleascollmeth:23} for $b=0$:
\begin{eqnarray*}
 	\mathcal{E}(u_{Z,\alpha}-u^*)
 	&\leq& C\Big( h_Z^{2m-4-\dS}  \big( \|\dot{u}^*\|_{  H^{{m-2}}(\M)} ^2 +\|u^*\|_{  H^m(\M)}^2 \big)
 	\\ && \qquad\qquad
 	+\big(  h_Z^{2m-2} + h_X^{-2}h_Z^{2m-\dS} +
 	h_X^{2m-2}  
 	\big)\| u^*(\cdot,0) \|_{H^m(\M)}^2
 	\Big)
 \end{eqnarray*}
with  the error functional defined as
\begin{equation}\label{eq:errorfunctional}
 	\mathcal{E}(u) :=
 	\esssup_{0\leq t\leq T}\|u\|_{H^1(\M)}^2
 	+\|u\|_{L^2(0,T;H^1(\M))}^2
 	+\|\dot u\|_{L^2(0,T;H^0(\M))}^2.
 \end{equation}
\green{Do note that the semi-discretized equation is not an ODE and the semi-discretized solution is not meant to be calculated.

The denseness requirement \eqref{eq:desnseness stability} serves to provide a sufficient (but not necessary) condition that ensures a stability estimate holds for all functions in the RBF trial space. The generic constant $C$ there  originates from those in a sampling inequality \cite[Lem. 3.1]{Cheung+LingETAL-leaskerncollmeth:18} and Bernstein's inverse inequality \cite[Thm. 10]{Hangelbroek+NarcowichETAL-DireInveResuBoun:17}. By explaining this, we highlight that the theory only guarantees the existence of such a $C$, without offering a practical way to estimate its value. Consequently, the concept of ``sufficiently dense'' (or overasampling ratio $n_X/n_Z$ for quasi-uniform  point sets) needs to be verified numerically.
}

\green{Next step in MoL is to define the fully discretized solution, which is defined by a difference equation, by discretizing in time using some $p$-th order backward difference scheme. We then show that the ODE, which is expected by everyone, is the limit of the difference equation. Formally, we redefine the semi-discretized solution}
in the form of \eref{eq:trial sp0} with coefficient functions   $\lambda_Z \in \big(C^{p+1}[0,T]\big)^{n_Z}$ \teal{satisfying} the $n_Z \times n_Z$ ODE system
\begin{equation}\label{eq: final ode}
 	\dot\lambda_Z(t) = \big[ - \Psi(X,Z)^\dagger\call_\M\Psi(X,Z)\big] \lambda_Z(t)   +   \Psi(X,Z)^\dagger   f(X,t),
 \end{equation}
where $\dagger$ denotes the pseudoinverse operator.  Alternatively, we can work with nodal values by solving
\begin{eqnarray}
 	\dot u_Z(Z,t) &=& \big[ - \Psi(Z,Z) \Psi(X,Z)^\dagger\call_\M\Psi(X,Z)  \Psi(Z,Z)^{-1} \big] u_Z(Z,t)
  \nonumber\\
   &&\qquad\qquad\qquad\qquad\qquad\qquad+   \Psi(Z,Z)  \Psi(X,Z)^\dagger   f(X,t),
   \label{eq: final ode2}
\end{eqnarray}
subject to an initial condition $\lambda_Z(0)$ \teal{or} $u_Z(0)$, respectively, that enjoys $H^1(\M)$-convergence. One possible approach is to use the interpolant of initial condition $g$ in \eref{eqPDE-b} from  the trial space $\calu_{Z,\M ,\Psi}$.
\green{
If the exact solution has $p+1$ continuous time derivatives, then the error estimate for the semi-discretized solution extends to the fully discretized solution with an additional temporal error term $\mathcal{O}(\triangle t^p)$, see the verification in our another work \cite[Tab. 1]{Chen+CheungETAL-kernleascollmeth:23}.
}


\green{
In theory, the mass matrix in the ODE is of full rank. However, numerically, zero eigenvalues may occur due to ill-conditioning. For the sake  of our discussion, we will still refer to these cases as ODEs, and this should not lead to any confusion.
Moreover,}
the ODE matrices in \eqref{eq: final ode} and \eqref{eq: final ode2} are similar matrices and, \teal{therefore}, share the same eigenvalues.  In our numerical studies, we only consider the former.
In our implementation, we do not use the ODE \eqref{eq: final ode} with a pseudoinverse matrix.
Instead, we compute the reduced QR-factorization of  the $n_X\times n_Z$ matrix $[\Psi( X, Z)] = : QR$, where $Q$ is $n_X\times n_Z$  and $R$ is square.
Then, we solve the equation
\begin{equation}\label{eqQRODE}
 	R \dot\lambda_Z(t) = \big[-Q^T\call_\M\Psi(X,Z)\big] \lambda_Z(t)+ Q^T f(X,t),
 \end{equation}
with the mass matrix $R$, which can be solved by many black-box ODE solvers and we use the explicit ODE45 Runge-Kutta method in MATLAB.
\green{
We choose  this integrator for its built-in adaptive time-stepping, which is a clear advantages of using MoL. It also allows us to infer the eigenvalue stability of the ODE. Additionally, the high order of convergence in time offered by the Runge-Kutta method  aligns well with the high order of convergence in space inherent to RBF.
}

We end this methodology section by remarking  that there are different ways to compute the entries in $[\call_\M\Psi(X,Z)]$:
\begin{enumerate}
 	\item If we know the normal vector $\n$ of $\M$ analytically, we can symbolically obtain $\call_\M\Psi$ using the formulas \eqref{def:grad} and \eqref{def:div}.
 	\item If $\n$ is known pointwise in some oriented point cloud of $\M$, we can use a pseudospectral RBF approximation \cite{Fuselier+Wright-HighKernMethDiff:13,Chen+Ling-Extrmeshcollmeth:20} to approximate $\call_\M\Psi$ at any point in the point cloud.
 	\item If $\M$ is implicitly defined by some discrete point sets, we can  first estimate all the necessary normal vectors using kernel interpolation level-set approaches \cite{Marchandise+PiretETAL-meshrepawithradi:12,Flyer+FornbergETAL-rolepolyRBF-appr:16}, and then return to Case 2.
\end{enumerate}
In our numerical experiments, we consider surfaces \teal{defined by} parametric \CMadd{or level-set} equations so that the normal vector $\n$ is known analytically.

\section{Oversampling vs Eigenvalue stability}\label{sec:Eigenvalue stability}
Over a decade ago, it was numerically shown that using least-squares in RBF interpolation \cite{Buhmann-Radibasifunc:03} and the method of lines can be beneficial \cite{Platte+Driscoll-Eigestabradibasi:06}. Recently, the convergence of RBF-LSC method for elliptic PDEs has been proven \cite{Cheung+LingETAL-leaskerncollmeth:18}. However,  oversampling is not yet a standard \teal{practice} in RBF method of lines, as it is not always necessary. For example, in certain square cases, \teal{the authors of \cite{2005Approximation,2012Solving,Dereli-meshkernmethline:12,Hussain+HaqETAL-Numesolusineequa:13} used square cases with $Z=X$ sucessfully}. Without oversampling, the DAE in \eqref{eq:DAE} become a system of ODEs with a symmetric positive definite mass matrix $\Psi(Z,Z)$. Therefore, while published results typically show successful cases, they may not provide a comprehensive understanding of the necessity of oversampling.

\green{The theoretical predicted convergence rate was verified in our previous work \cite[Exmp. 1 \& 2]{Chen+CheungETAL-kernleascollmeth:23} and will not be covered again here. }
The rest of this paper presents numerical experiments aimed at identifying an appropriate range, if needed at all, \teal{for} oversampling ratio $n_X/n_Z$ that yields stable eigenvalues. 
\green{We will make numerical observations aiming at establishing a rule-of-thumb for effectively configuring the RBF-LSC-MoL for practical use.}
In all examples below, we use the $H^{m+1/2}(\R^d)$-reproducing Sobolev kernel  with on-surface smoothness order $m>(d-1)/2$ defined by
\[
 \Psi_{m+1/2}(x,z) := \|x-z\|_{\ell^2(\R^d)}^{(m+1/2)-d/2}\calk_{(m+1/2)-d/2}\big( \|x-z\|_{\ell^2(\R^d)} \big)\quad \text{for $x,z\in\R^d$}
 \]
using the modified Bessel functions of the second kind $\calk$.

\exmp{Isotropic diffusion-reaction on the unit sphere}{ex1}
\green{The focus of this exploration will be to understand the intricate interplay between various numerical configurations. In particular, we will examine the effects of oversampling ratio, as well as the influence of the smoothness orders $m$ of the kernels and numbers of trial centers.}
 \begin{figure}
 	\centering
 	\includegraphics[width=0.5\textwidth]{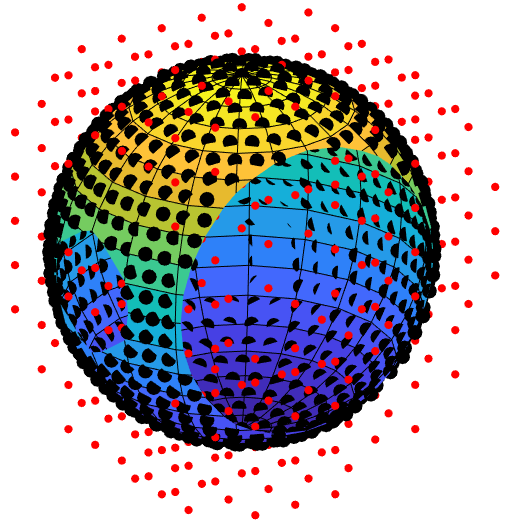}
 	\caption{Exmp.~\ref{ex1}--\ref{ex2}: Schematic demonstration of data point distributions. (Black) Quasi-uniform data points on the unit sphere, and (Red) Regular grid in some narrow band domain containing the unit sphere.}\label{fig0}
 \end{figure}

We consider the surface PDE $\dot u -  \Delta_\M  {u} + 3 u = f$ on the unit sphere. We generate \cite{Ling-poin:14} sets
of quasi-uniform RBF centers $Z\subset\M$ of size $n_Z$, where
 \begin{equation}\label{nz}
 n_Z\in\{658, 1266, 2506, 3850\},
 \end{equation}
\teal{as shown by} the on-surface black dots in Figure~\ref{fig0} for a schematic demonstration.
The irregular pattern in the tested $n_Z$ is for the sake of comparison in the next example.
We test oversampling ratios
\[
n_X/n_Z\in\{1,1.5,2,4\},
\]
and generate sets of $n_X$ collocation points for each tested set $Z$ in a similar manner.
We compute all eigenvalues of the ODE matrix $\big[ - \Psi(X,Z)^\dagger\call_\M\Psi(X,Z)\big]$ in \eqref{eq: final ode}
using the constructed sets of $X$ and $Z$. We then present a summary of our observations to the readers.

 \begin{figure}
 	\centering
 	\begin{overpic}[width=0.9\textwidth]{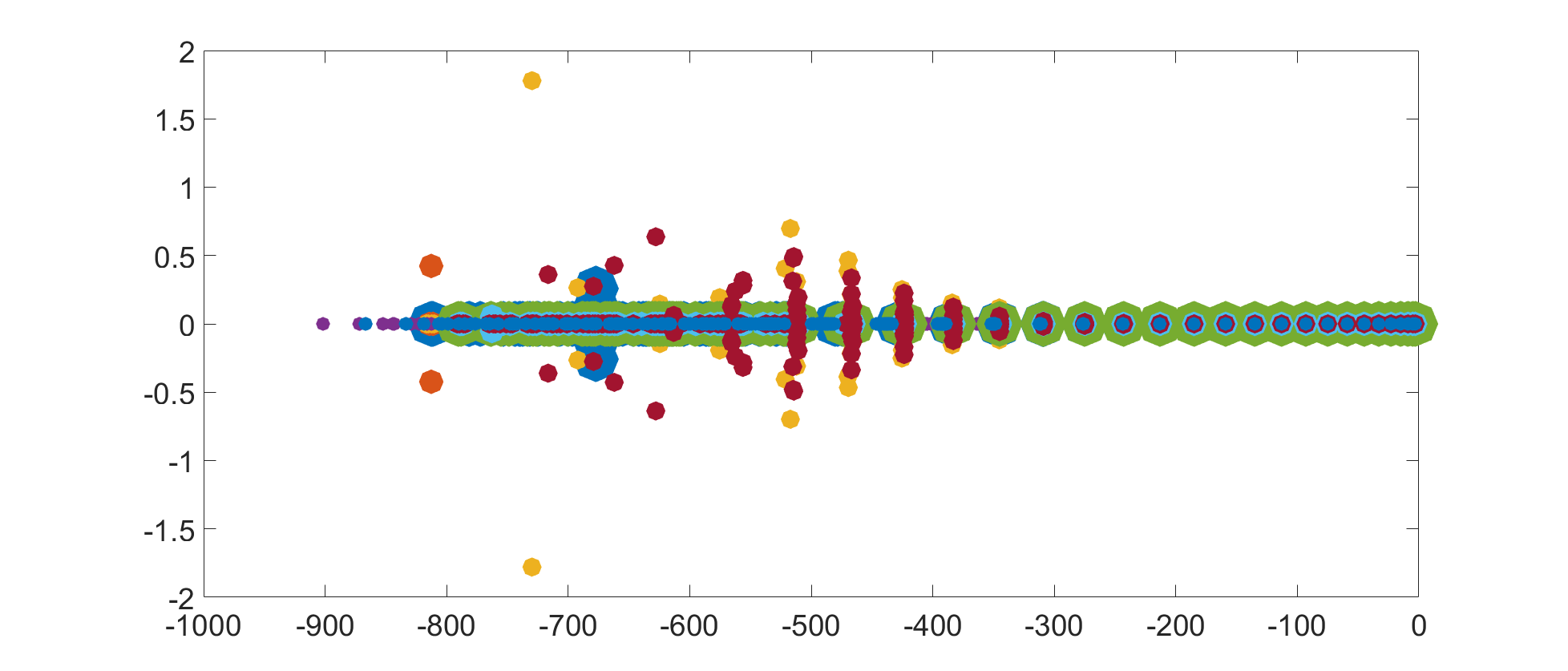}
 		\put(6,20){\scriptsize \rotatebox{90}{$\mathfrak{Im}$}}
 	\end{overpic}
 	\\
 	\begin{overpic}[width=0.9\textwidth]{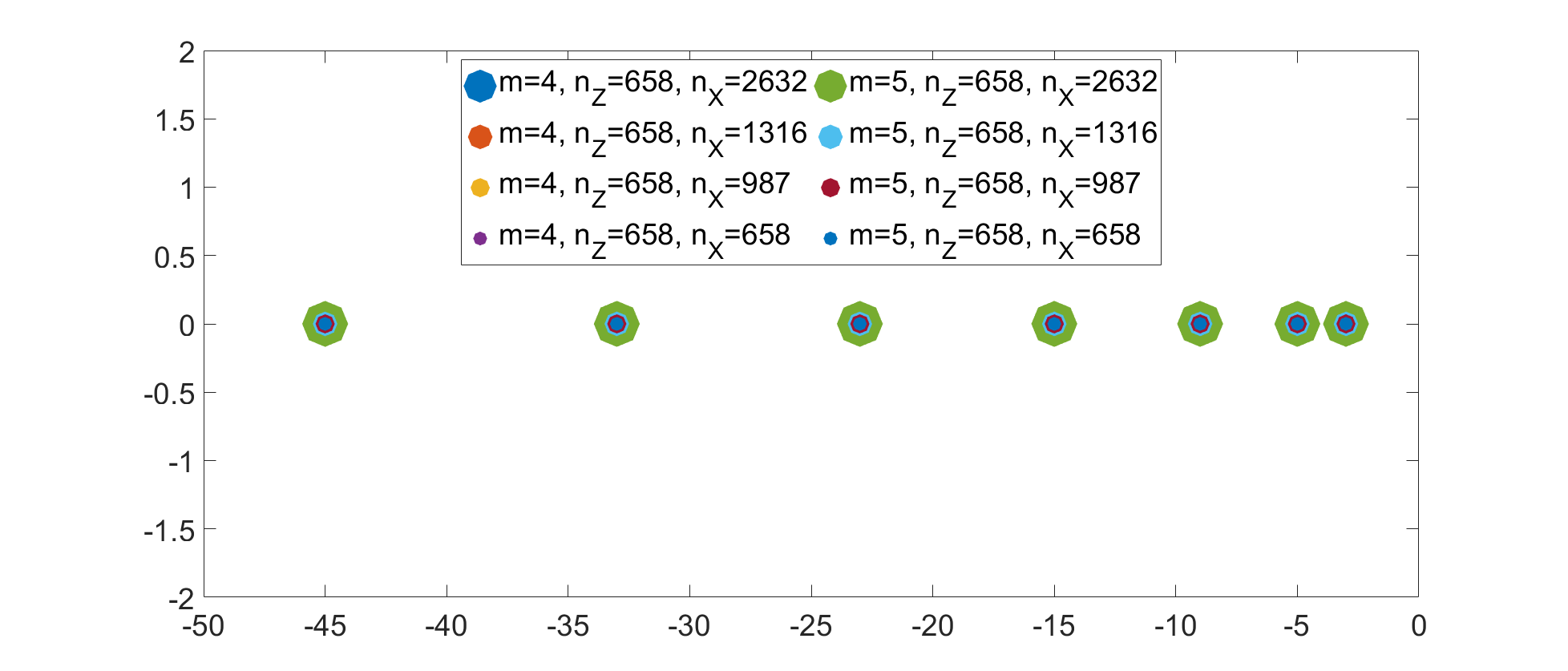}
 		\put(50,42){\scriptsize \rotatebox{0}{$\mathfrak{Re}$}}
 		\put(6,20){\scriptsize \rotatebox{90}{$\mathfrak{Im}$}}
 		\put(50,0){\scriptsize \rotatebox{0}{$\mathfrak{Re}$}}
 	\end{overpic}
 	\caption{Exmp. \ref{ex1}: We observe stable eigenvalues in $\C$ when both the kernel order of smoothness $m$ and the number of RBF centers $n_Z$ are relative small; oversampling does not have any obvious effect in terms of eigenvalue distribution. The provided legend applies to both subfigures: (TOP) full spectra, and (BOTTOM) Zoom-in near the origin. }\label{fig1}
 \end{figure}

 \subsection*{Observation 1: {Oversampling is not needed in cases of small $m$.}}

 In Figure~\ref{fig1}, we present the eigenvalue distributions corresponding to the ODE matrices with the smallest $n_Z=658$ RBF centers and smoothness orders $m=4$ and $5$, which were selected due to their stable eigenvalues. In the top figure, we observe that a higher oversampling ratio (indicated by larger markers in the figure) can reduce the spectral radius and   result  in eigenvalues closer to the real number line. Zooming in towards the origin in the bottom figure, we can see that all eight test cases can  approximate the largest eigenvalues well.

 \begin{figure}
 	\centering
 	\begin{overpic}[width=0.9\textwidth]{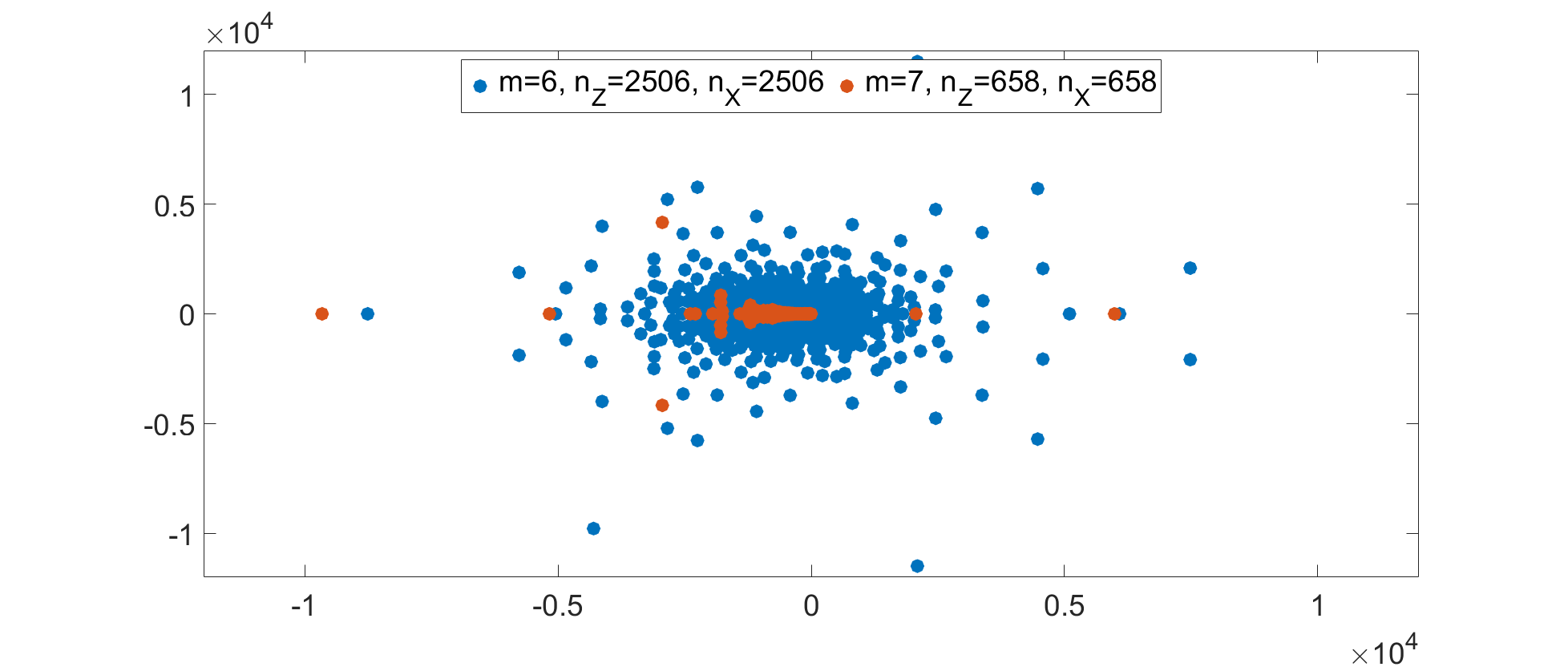}
 		\put(6,20){\scriptsize \rotatebox{90}{$\mathfrak{Im}$}}
 		\put(50,0){\scriptsize \rotatebox{0}{$\mathfrak{Re}$}}
 	\end{overpic}
 	\caption{Exmp. \ref{ex1}: RBF method of lines without any oversampling has unstable eigenvalues when either $m$ and/or $n_Z$ increases.}\label{fig2}
 \end{figure}

 \subsection*{Observation 2: {Eigenvalues become unstable if either $m$ or $n_Z$ increases.}}
An interesting observation, \teal{for} which we have not provided evidence, is that for $m=6$, eigenvalues \teal{remain} stable for both $n_X=n_Z=658$ and $1266$ without oversampling. However, when we keep $m=6$ fixed and increase the number of RBF centers to $n_Z=2506$, instability arises in the eigenvalue distributions, as \teal{indicated} by the blue markers in Figure~\ref{fig2}. Similarly, if we keep $n_Z=658$ fixed \teal{as shown} in Figure~\ref{fig1} but increase the smoothness order to $m=7$, we \teal{once again} observe unstable eigenvalues, indicated by orange markers in the figure. \teal{This suggests} that the \teal{RBF-LSC-MoL} does not allow for $m$- or $n_Z$-convergence without oversampling.


 \begin{figure}
 	\centering
 	\begin{overpic}[width=0.9\textwidth]{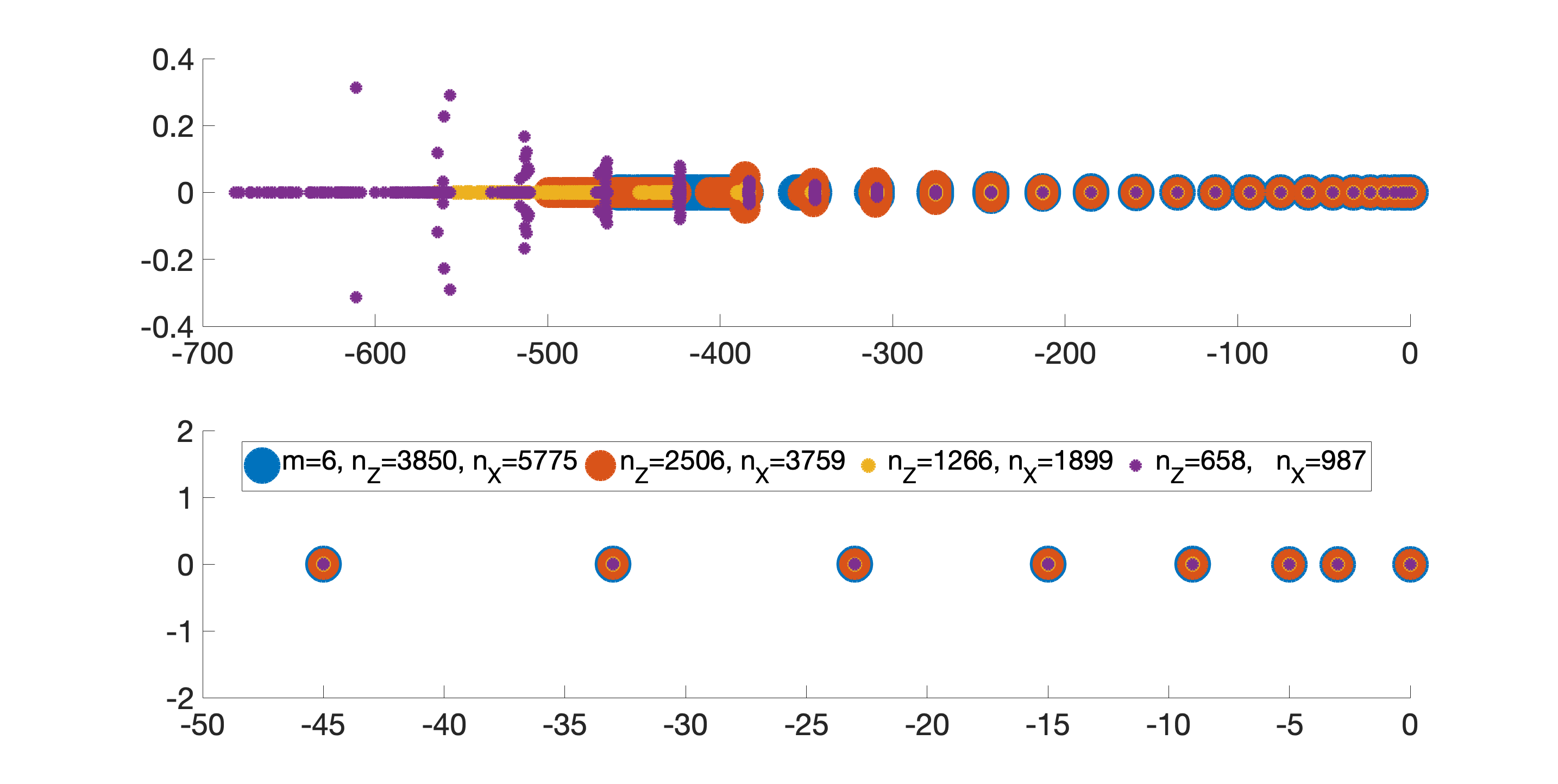}
    \put(2,25){ \rotatebox{90}{ {$m=6$}  }}
 	\end{overpic}

\bigskip
	\begin{overpic}[width=0.9\textwidth]{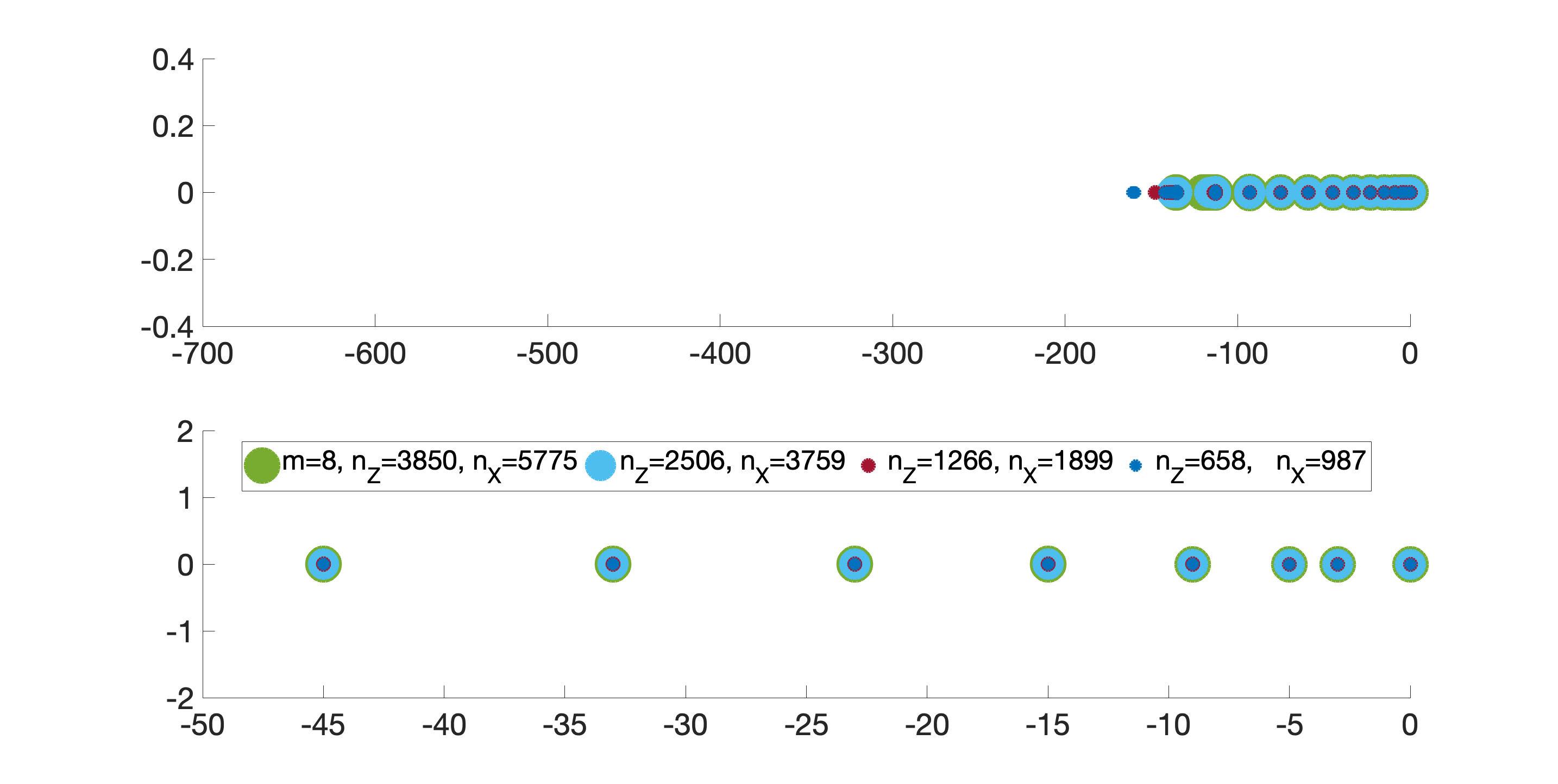}
\put(2,25){ \rotatebox{90}{ {$m=8$}  }}
\end{overpic}

\bigskip
 \begin{overpic}[width=0.9\textwidth]{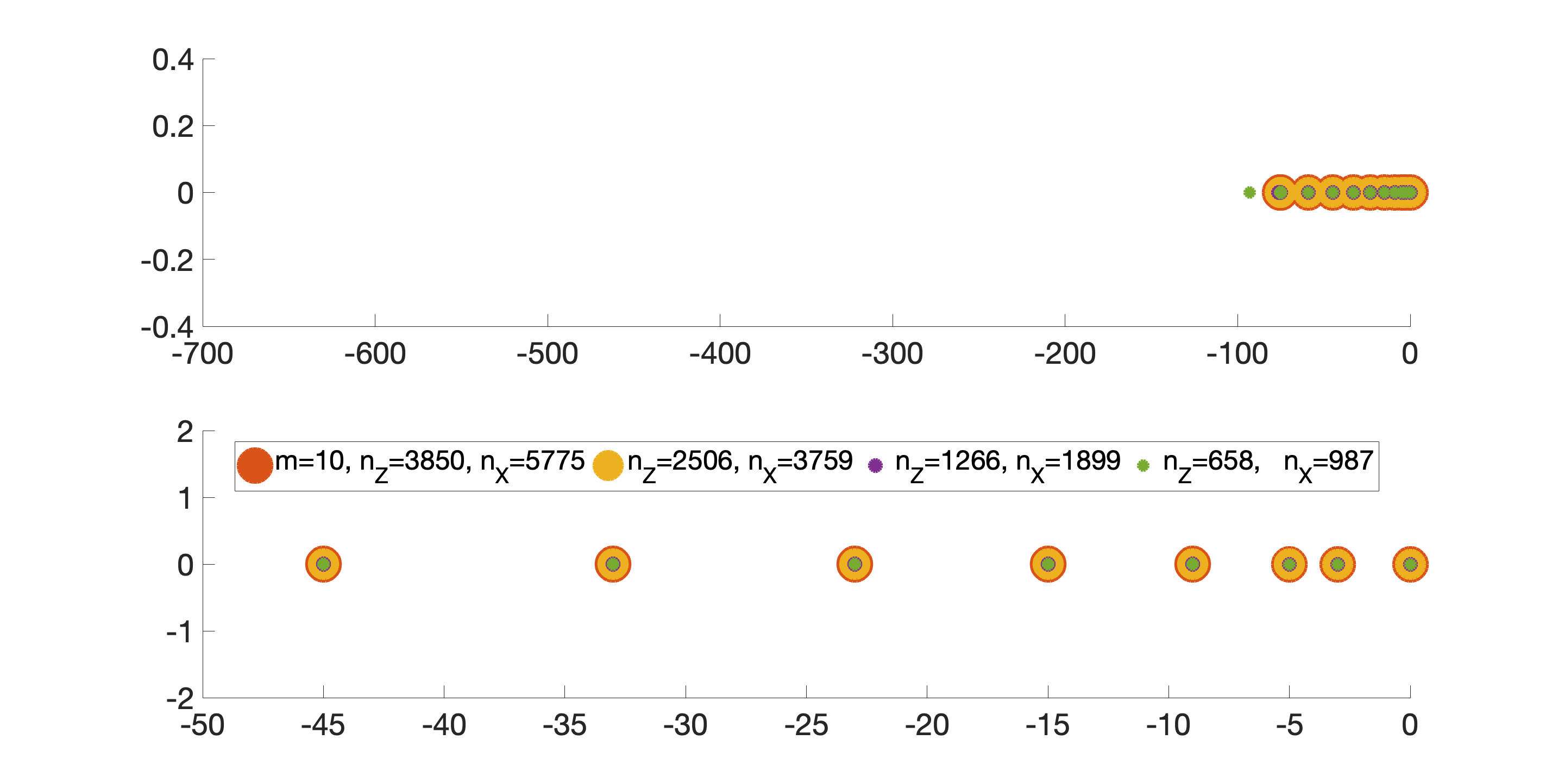}
\put(2,25){ \rotatebox{90}{ {$m=10$}  }}
\end{overpic}

 	\caption{ Exmp. \ref{ex1}: By imposing $1.5\times$ oversampling ratio,  we consistently observe stable eigenvalues in $\C$ when using kernel with high order of smoothness and more RBF centers. Do  note the spurious zero eigenvalue. 
 \teal{Subfigures for $m=6$, $8$, and $10$ here are presented in same style as Figure~\ref{fig1}.}
 }\label{fig3}
 \end{figure}

 \subsection*{Observation 3:  {Stability in eigenvalues can be regained by oversampling.}}
For all tested $n_Z$, we \teal{consider cases where the oversampling ratio is greater than one, and avoid the $X=Z$ case.}
Figure~\ref{fig3} shows the eigenvalues obtained with $m=6,8,10$, while $m=7,9$ were omitted due to their high similarity.
By comparing these results with the top of Figure~\ref{fig1}, we see that the spectra in this example lie in a smaller range. Yet, the most significant and noteworthy difference is a spurious zero eigenvalue that is clearly visible in the bottom subfigure.
We report that all 12 tested cases in this example contain a zero eigenvalue with different multiplicity, depending on $m$. This is  {the problem of ill-conditioning}  in RBF methods when we use the RBF directly to define basis functions. The kernel matrix $[\Psi(X,Z)]$ eventually becomes \green{numerically} rank deficient as we increase the order $m$ and/or the number of points in $Z$. Thus, as we push the RBF approximation power to its finite-precision limit, the ODE in \eqref{eqQRODE} unavoidably becomes \green{degenerate}. The inverse relationship between condition number and accuracy, a.k.a. uncertainty principle \cite{Schaback-Erroesticondnumb:95}, for interpolation and the folklore for time-independent PDEs both suggest that we should tune the RBFs to yield nearly singular kernel/PDE matrices under the given computational precision. This raises the question of whether we should apply a similar approach to \teal{RBF-LSC-MoL}.


 \begin{figure}
 	\centering
 	\begin{overpic}[width=0.9\textwidth]{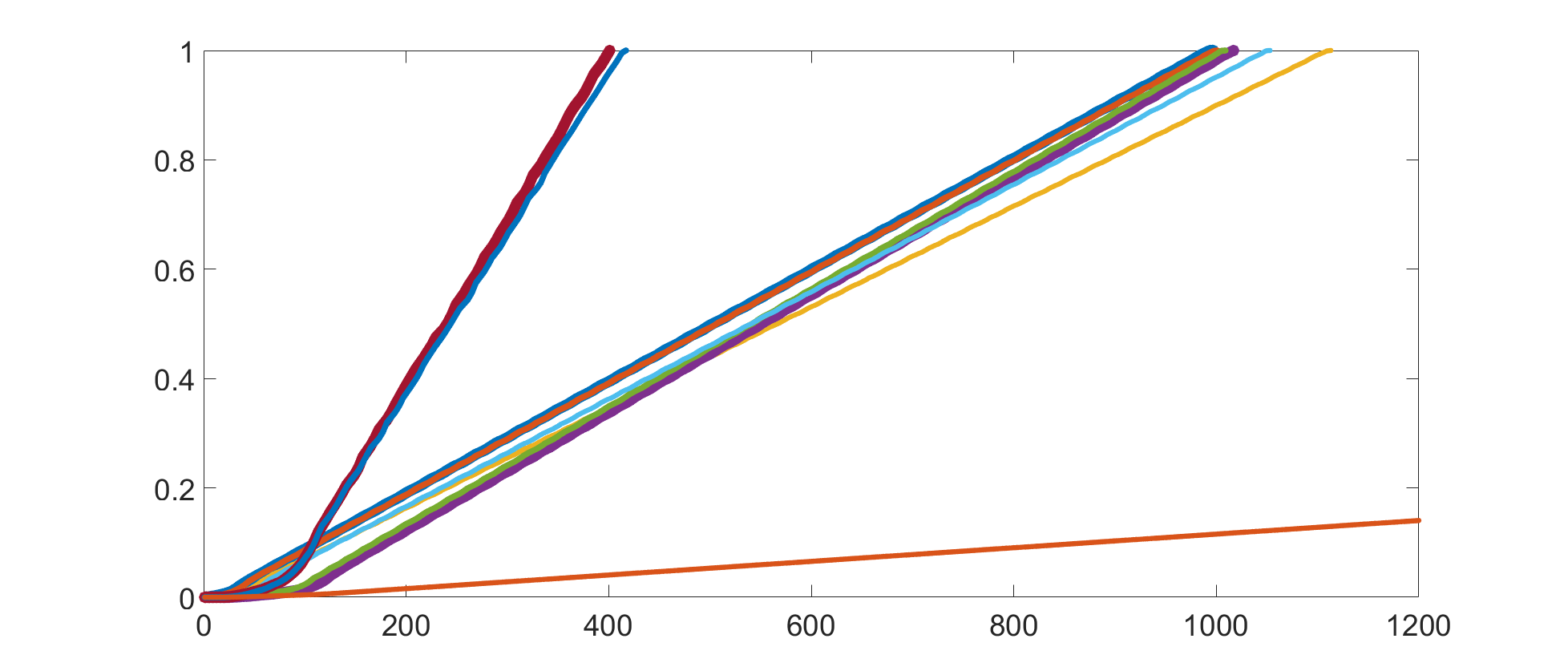}
 		\put(6,18){\scriptsize \rotatebox{90}{Time $t$}}
 		\put(44,1){\scriptsize \rotatebox{0}{\#\,Step}}
 	\end{overpic}
 	\\
 	\begin{overpic}[width=0.9\textwidth]{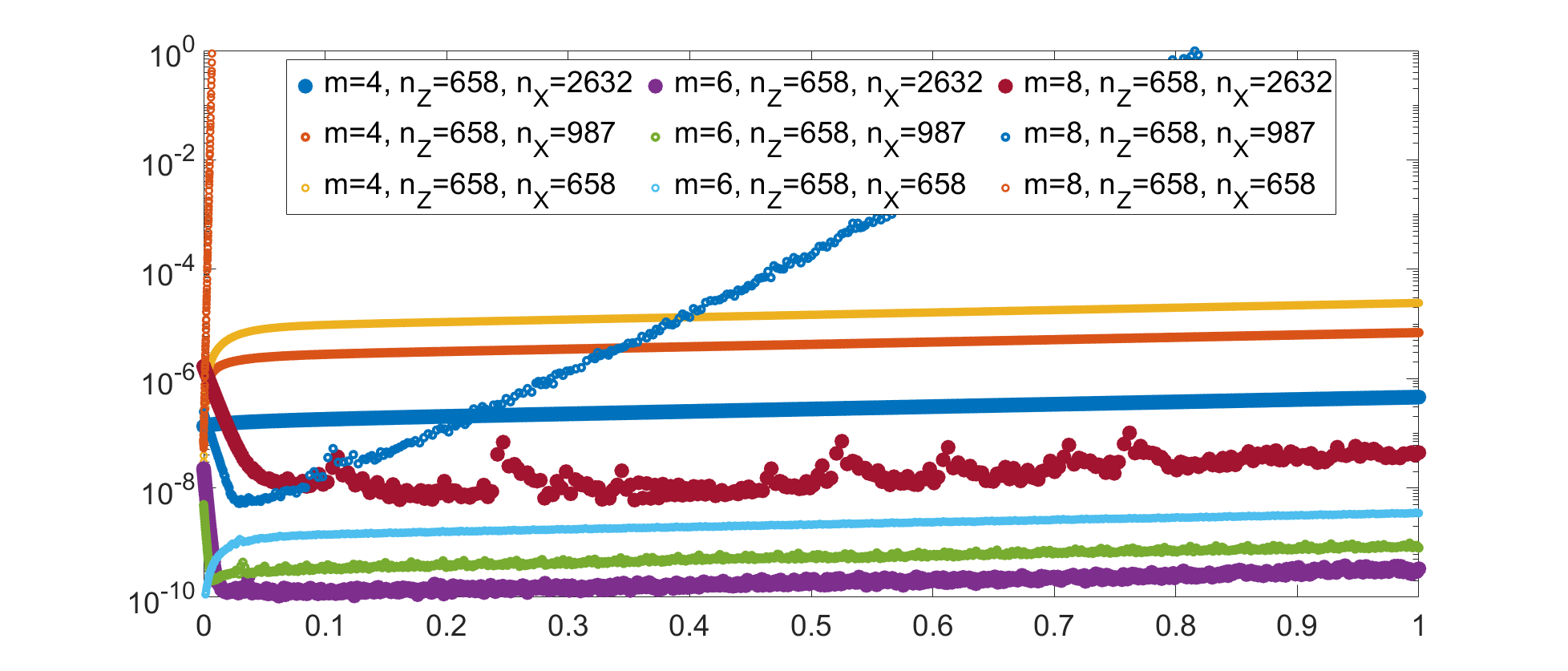}
 		\put(6,18){\scriptsize \rotatebox{90}{$L^\infty(\M)$--error}}
 		\put(44,0){\scriptsize \rotatebox{0}{Time $t$}}
 	\end{overpic}
 	\caption{Exmp. \ref{ex1}: Summaries of ODE45 results.
 		(TOP) Current time against the number of adaptive time steps,
 		(BOTTOM) maximum errors of various RBF method of lines solutions against time.
 	}\label{fig4}
 \end{figure}

 \subsection*{Observation 4: {\green{Degenerated ODE does not pose a problem to accuracy}.}}

\teal{To compute accuracy, we pick}  $u^*(x,t) =\exp( x_1+1/(1+t) )$ for $x=(x_1,x_2,x_3)\in\M$ and  $t\in[0,1]$ be the solution to our target PDE. We compute the corresponding source function $f := \partial_t u^* -  \Delta_\M  {u^*} + 3 u^* $  analytically.
 We then set up the ODE \eqref{eqQRODE} using a fixed $Z$ with $n_Z=658$, but varying the kernel smoothness and oversampling ratios.
 All ODEs, some of which are \green{degenerate} due to the previous observation, are solved using ODE45 in MATLAB without fine-tuning any error tolerances. Representative results are presented in Figure~\ref{fig4}, which we will walk through in detail.

The top subfigure of Figure~\ref{fig4} shows the cumulative time taken by the time-stepping in ODE45. Note that, due to the adaptivity of the time-stepping, some of these curves are not straight lines, especially for the earlier steps. The two steepest curves, corresponding to $m=8$ and $n_X=987$ and $2632$, require the fewest time steps to complete the time integration and are the most efficient in terms of computation time. Since all $n_Z\times n_Z$ ODE matrices in this test are of the same size, requiring fewer time steps can translate into a faster overall PDE solving time. The extra cost of using different $n_X>n_Z$ is due to the computation of the reduced QR in \eqref{eqQRODE}, which grows linearly and quadratically with respect to $n_X$ and $n_Z$, respectively \cite[Tab. 3.1]{Ling-fastblocalgoquas:16}. The outlier on the other end of the spectrum also comes from $m=8$ but without oversampling. This setup has unstable eigenvalues, causing ODE45 to run with very small time steps in an attempt to solve this stiff problem  and numerical garbage will be produced  after a long runtime. The group of curves in the middle corresponds to results obtained using $m=4$ and $6$ with all tested oversampling ratios.


 \begin{table}[]
 	\centering
 	\begin{tabular}{@{}l|ll|ll|ll|ll|ll|@{}}
 	\toprule
 		& \multicolumn{2}{|c|}{$m=4$}                           & \multicolumn{2}{c|}{$m=5$}                           & \multicolumn{2}{c|}{$m=6$}                           & \multicolumn{2}{c|}{$m=7$}                           & \multicolumn{2}{c|}{$m=8$}                            \\ \hline
 		\multicolumn{1}{|l|}{$n_X$} & \multicolumn{1}{l|}{$\kappa$} & \multicolumn{1}{l|}{Err} & \multicolumn{1}{l|}{$\kappa$} & \multicolumn{1}{l|}{Err} & \multicolumn{1}{l|}{$\kappa$} & \multicolumn{1}{l|}{Err} & \multicolumn{1}{l|}{$\kappa$} & \multicolumn{1}{l|}{Err} & \multicolumn{1}{l|}{$\kappa$} & \multicolumn{1}{l|}{Err.} \\ \hline
 		\multicolumn{1}{|l|}{658}   & 5e9  & 2.4e-5  & 1e12  & 3.0e-7  & 6e14   & 3.4e-9   & 6e18   & NaN      & 6e19   & NaN     \\
 		\midrule
 		\multicolumn{1}{|l|}{987}   & 8e9  & 6.9e-6  & 2e12  & 8.2e-8  & 8e14   & 8.0e-10  & 1e17   & 7.6e-10  & 2e17   & 8.1e+1  \\
 		\midrule
 		\multicolumn{1}{|l|}{\teal{1316}}   & 8e9  & 9.5e-7  & 2e12  & 1.3e-8  & 8e14   & 3.1e-10  & 1e17   & 8.0e-10  & 1e17   & 2.4e-4  \\ \midrule
 		\multicolumn{1}{|l|}{2632}  & 8e9  & 4.6e-7  & 2e12  & 1.1e-8  & 8e14   & 3.3e-10  & 9e16   & 8.1e-10  & 1e17   & 4.4e-8  \\ \botrule
 	\end{tabular}
 	\caption{Exmp. \ref{ex1}: With $n_Z=658$ fixed, the condition number ($\kappa$) of the kernel matrix $\Psi(X,Z)$ and final $L^\infty(\M)$-error of the RBF method of lines solution at the final time $t=1$.}\label{tab0eval}
 \end{table}

We now turn to the bottom of Figure~\ref{fig4} to study the accuracy of numerical solutions obtained from the RBF-LSC-MoL. From smallest to largest error at $t=0.1$, we observe that:
 \begin{itemize}
 	\item three $m=6$ error curves (in colors purple, green, and sky blue) have the smallest errors,
 	\item two $m=8$ error curves (in colors red and cobalt) comes next,
 	\item $m=8$ case  without oversampling blows the errors up immediately.
 \end{itemize}
We also consider the condition number of the kernel matrix $\Psi(X,Z)$ instead of that of the ODE matrix in \eqref{eq: final ode}. Note that the condition number of the ODE matrix is the square of its kernel matrix counterpart, which is listed in Table~\ref{tab0eval}.
\green{When $m=4$, there are no (numerically) zero eigenvalues in the kernel matrix $\Psi(X,Z)$ {, as shown in} the bottom plot of Figure~\ref{fig1}. In this example, we are dealing with full rank {ODE  matrices}}. On the other hand, when $m=6$ and $8$, we are working with \green{degenerate ODEs} subject to different degrees of rank loss in the ODE matrix. Recall the eigenvalue distribution for these cases in Figure~\ref{fig3}.



We observe that the least error in Table~\ref{tab0eval} occurs when $\kappa(\Psi(X,Z)) \lesssim 1/\varepsilon_{\text{machine}}$. Thus, the folklore  remains helpful in selecting the smoothness order $m$. Yet, indicators for selecting $n_X$ are still missing. From the last column ($m=8$) of the table, we see that the solution accuracy heavily depends on the oversampling ratio. However, the ranks and condition numbers of all three ODE matrices $\big[ - \Psi(X,Z)^\dagger\call_\M\Psi(X,Z)\big]$ with $987\leq n_X \leq 2632$ are 121 and $\calo(\text{1e17})$, respectively.

\bigskip
\exmp{Using regular RBF centers $Z$ in some narrow band domain $\Omega\supset\M$}{ex2}
\green{In this second example, we venture beyond the theoretical constraints regarding the distribution of trial centers.}
Given the same number $n_Z$ of RBF centers, we can maximize the minimum separating distance $q_Z$ of $Z$ if we place $Z$ out-of-surface. However, we emphasize that our convergence theory in \cite{Chen+CheungETAL-kernleascollmeth:23} does not apply in this case, although this is the data point structure used in embedding methods \cite{Ruuth+Merriman-SimpEmbeMethSolv:08,Piret-orthgradmeth:12,Cheung+Ling-Kernembemethconv:18}.

 \begin{figure}
 	\centering
 	\begin{overpic}[width=0.9\textwidth]{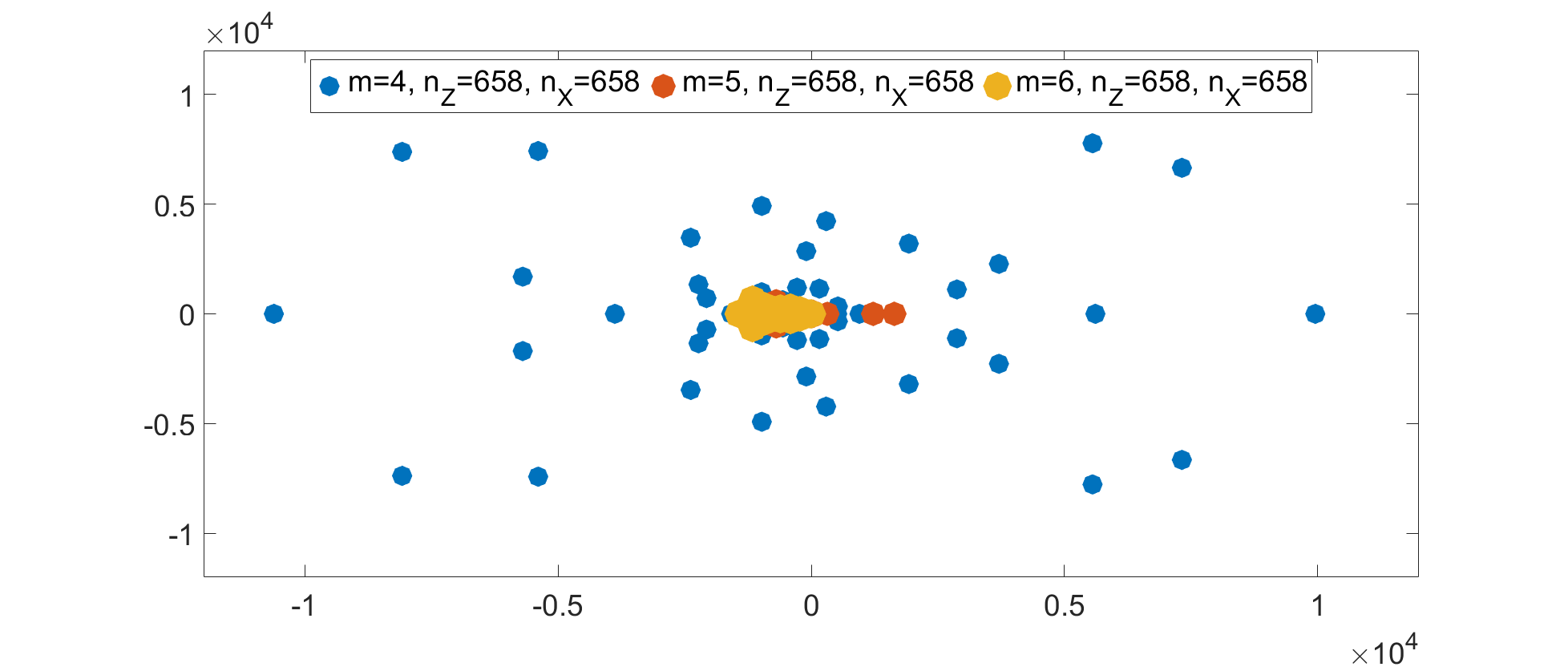}
 		\put(6,20){\scriptsize \rotatebox{90}{$\mathfrak{Im}$}}
 	\end{overpic}
 	\\
 	\begin{overpic}[width=0.9\textwidth]{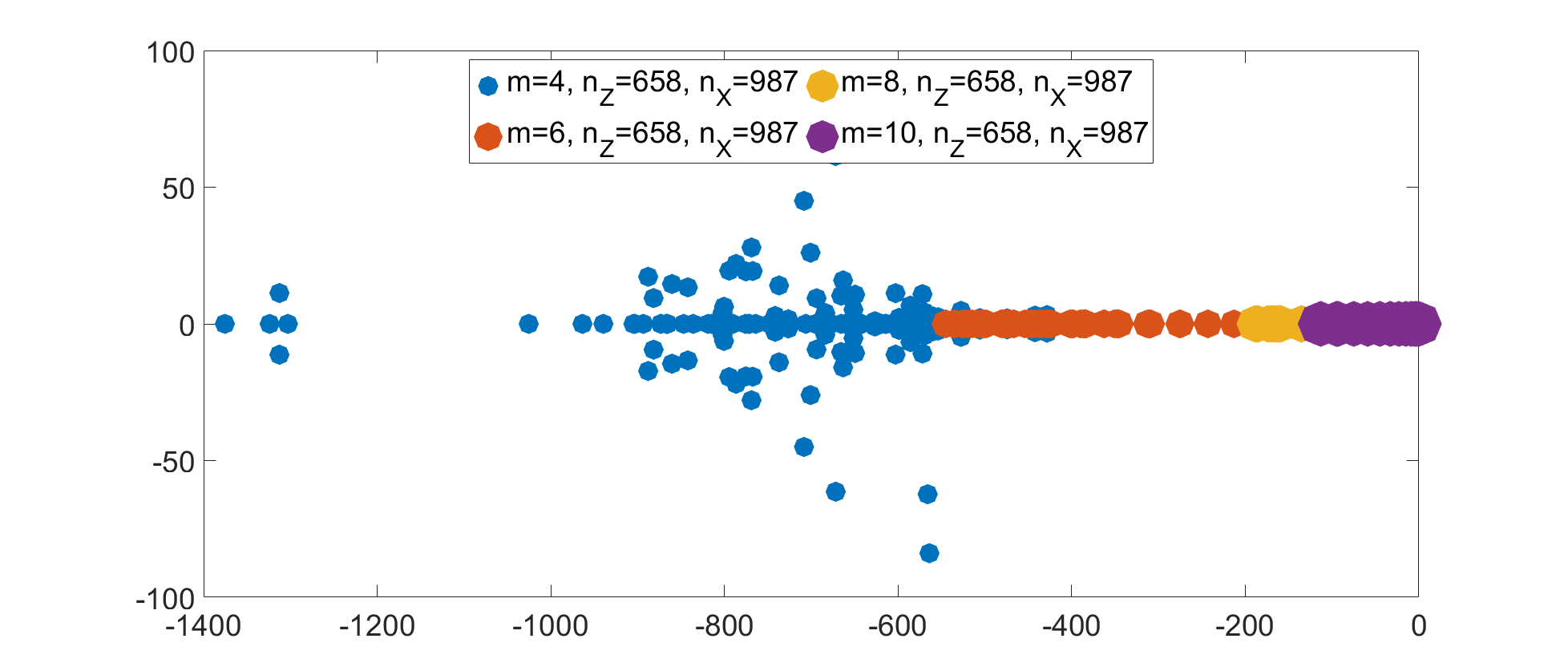}
 		\put(50,42){\scriptsize \rotatebox{0}{$\mathfrak{Re}$}} \put(6,20){\scriptsize \rotatebox{90}{$\mathfrak{Im}$}}
 		\put(50,0){\scriptsize \rotatebox{0}{$\mathfrak{Re}$}}
 	\end{overpic}
 	\caption{Exmp. \ref{ex2}: When we place RBF centers in a narrow band domain containing $\M$, (TOP) square approach with $X=Z$ yields unstable eigenvalues for all tested smoothness orders $m$, and (BOTTOM) imposing $1.5\times$ oversampling ratio is sufficient for stable eigenvalues.}\label{fig5}
 \end{figure}

For some small parameter $\delta$, let the narrow band domain be
\[
 \Omega := \Omega_{\delta,\M}= \Big\{ x\in\R^d:\, \|x-z\|_2\leq \delta \text{ for some $z\in\M$}\Big\}.
 \]
We define the set of near-surface regular RBF centers as $Z:= (\green{(1\pm\delta)}\mathbb{W} )^d \cap \Omega$. By using four values of $\delta$ between 0.125 and 0.25, we can obtain sets of $Z$ with exactly $n_Z$ points as in \eqref{nz} in Example~\ref{ex1}, showed as red dots in Figure~\ref{fig0}. Collocation points $X$ must lie on the surface, and they are generated as in the previous example.

In Figure~\ref{fig5}, we see from the top subfigure that unstable eigenvalues of the ODE matrices appear in all tested values of $m$, including $m=4,5$, which worked fine in the previous example with $Z\subset\M$, when oversampling is not used. In the bottom subfigure, oversampling can greatly stabilize the eigenvalue distribution.

Turning to Table~\ref{tab1eval}, we note that $m=6$ is still the optimal order of kernel smoothness. We also observe that the error profiles for $m=7,8$ behave very similarly to $m=8$ in the previous example. In all cases, the $n_X$-convergence behavior is pronounced.

 \begin{table}[]
 	\begin{tabular}{@{}l|ll|ll|ll|ll|ll|@{}}
 	\hline
\multicolumn{1}{|l|}{ }  & \multicolumn{2}{c|}{$m=4$}                           & \multicolumn{2}{c|}{$m=5$}                           & \multicolumn{2}{c|}{$m=6$}                           & \multicolumn{2}{c|}{$m=7$}                           & \multicolumn{2}{c|}{$m=8$}                            \\ \hline
 		\multicolumn{1}{|l|}{$n_X$} & \multicolumn{1}{l|}{$\kappa$} & \multicolumn{1}{l|}{Err} & \multicolumn{1}{l|}{$\kappa$} & \multicolumn{1}{l|}{Err} & \multicolumn{1}{l|}{$\kappa$} & \multicolumn{1}{l|}{Err} & \multicolumn{1}{l|}{$\kappa$} & \multicolumn{1}{l|}{Err} & \multicolumn{1}{l|}{$\kappa$} & \multicolumn{1}{l|}{Err.} \\ \hline
 		\multicolumn{1}{|l|}{658}   & 2e19   & NaN     & 3e19  & NaN     & 1e19  & NaN     & 3e19  & NaN     & 4e19  & NaN     \\ \hline
 		\multicolumn{1}{|l|}{987}   & 2e12   & 2.6e-6  & 3e14  & 1.5e-7  & 1e17  & 8.3e-9  & 2e17  & 9.4e-1  & 3e17  & 1.4e-1  \\ \hline
 		\multicolumn{1}{|l|}{1316}   & 1e12   & 2.4e-6  & 3e14  & 1.5e-7  & 7e16  & 7.2e-9  & 2e17  & 6.3e-5  & 2e17  & 1.1e-5  \\ \hline
 		\multicolumn{1}{|l|}{2632}  & 1e12   & 2.5e-6  & 3e14  & 1.5e-7  & 6e16  & 7.4e-9  & 1e17  & 7.8e-8  & 1e17  & 1.2e-8  \\ \hline
 	\end{tabular}
 	\caption{Exmp.~\ref{ex2}: Results  corresponding to Table~\ref{tab0eval} for near-surface regular RBF center in a narrow band domain.}\label{tab1eval}
 \end{table}

\bigskip
\exmp{Anisotropic diffusion with non-identity diffusion tensor}{ex3}
 \begin{figure}
 	\centering
 	\includegraphics[width=0.5\textwidth]{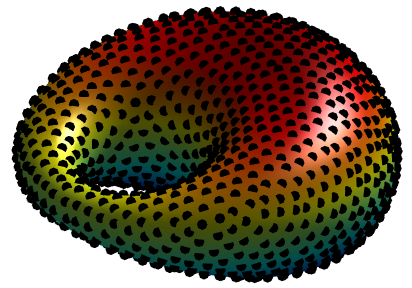}
 	\caption{Exmp.~\ref{ex3}: Schematic demonstration of quasi-uniform data point distributions on a Dupin Cyclide.}\label{figDC point}
 \end{figure}
\green{
Our theoretical framework is general enough to accommodate the Laplacian in divergence form, whose diffusion tensor should satisfy \cite[Assumption 1]{Chen+CheungETAL-kernleascollmeth:23}. In this third example, we aim to provide robust numerical evidence by working on a PDE with an anisotropic, i.e., space dependent, diffusion tensor}. For PDEs, we consider the differential operator \eqref{equLs} with reaction coefficient $b=3$ and a non-identity time-independent diffusion tensor
 \[
 A(x,t):= P(x)\,\text{diag}\big([1+x_1^2,\,1,\,1]\big),
 \]
where $P$ is the projection matrix in \eqref{def:P(y)}, so that the matrix $A(x,t):\TM\to\TM$ is symmetric positive definite on the tangent space $\TM$ of $\M$.
\teal{For centers distribution,} we considered a more complicated setup. Firstly, we place $n_Z=314,\,744,\,1320$ RBF centers and $n_X=2976,\,5296$ collocation points quasi-uniformly on a dupin cyclide, as shown in Figure~\ref{figDC point}. We only consider Sobolev kernels with smoothness $m=6$.

 \begin{figure}[]
 	\begin{tabular}{cc}
 		\begin{overpic}[width=0.48\textwidth,trim=48 48 48 48]{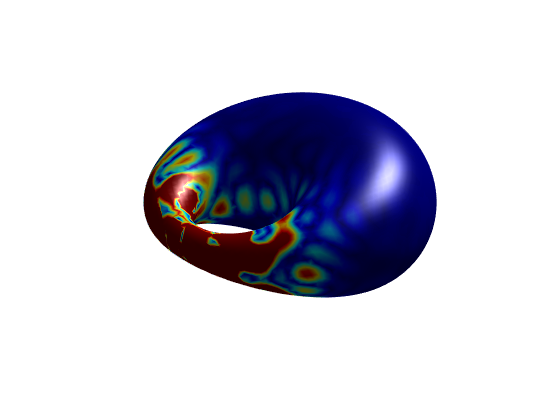}
 			\put(05,30){\rotatebox{90}{$n_Z=314$}}
 			\put(40,70){\rotatebox{0}{$n_X=2976$}}
 		\end{overpic}
 		&
 		\begin{overpic}[width=0.48\textwidth,trim=48 48 48 48]{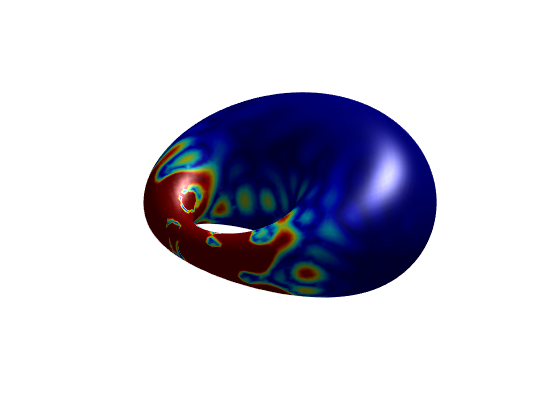}
 			\put(40,70){\rotatebox{0}{$n_X=5296$}}
 		\end{overpic}
 		\\
 		\begin{overpic}[width=0.48\textwidth,trim=48 48 48 48]{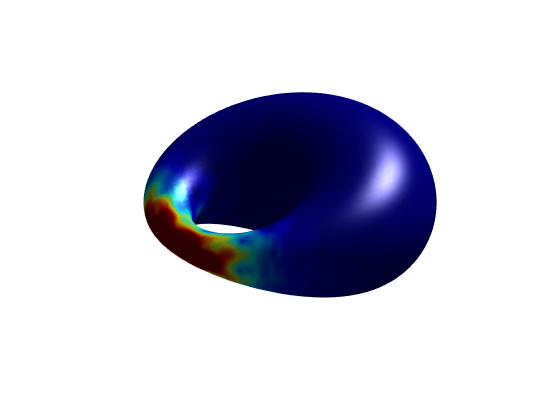}
 			\put(05,30){\rotatebox{90}{$n_Z=744$}}
 		\end{overpic}
 		&
 		\includegraphics[width=0.48\textwidth,trim=48 48 48 48]{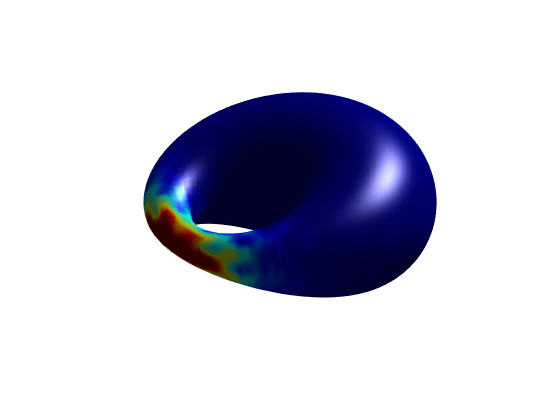}             \\
 		\begin{overpic}[width=0.48\textwidth,trim=48 48 48 48]{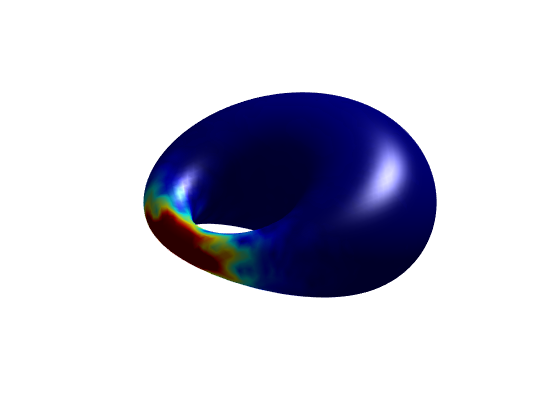}
 			\put(05,30){\rotatebox{90}{$n_Z=1320$}}
 		\end{overpic}
 		&
 		\includegraphics[width=0.48\textwidth,trim=48 48 48 48]{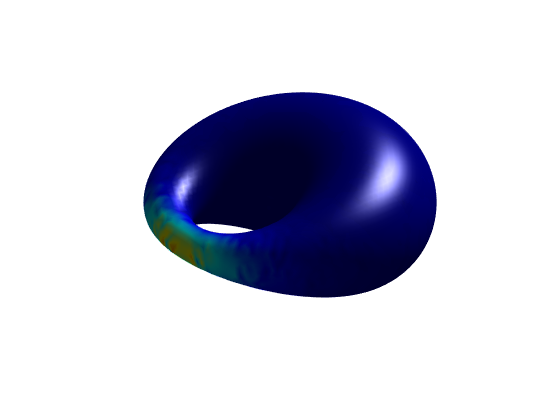}             \\
 		\multicolumn{2}{c}{\includegraphics[width=0.8\textwidth]{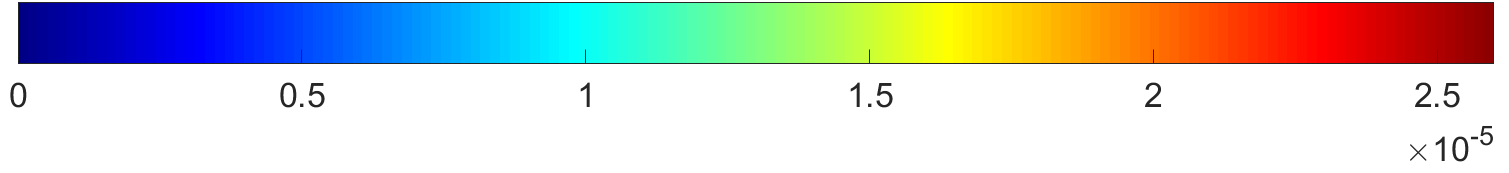}}
 	\end{tabular}
 	\caption{Exmp.~\ref{ex3}: Error profiles, which use the same colormap for easy comparison, of the RBF method of lines solutions to an anisotropic diffusion reaction equation for various numbers of RBF centers $n_Z$ (down the rows) and collocation points $n_X$ (across the columns). }\label{figDC Err}
 \end{figure}

We show all error profiles at time $t=1$ in Figure~\ref{figDC Err}. For a very small $n_Z=314$ of RBF centers, we see in the first row that 9.4$\times$ and 16.8$\times$ oversamplings yield no improvement in accuracy. Indeed, the error profile of $(n_Z,n_X) = (314,744)$ (omitted in the figure) also \teal{exhibits} a similar error pattern.
If we read Figure~\ref{figDC Err} \teal{from top to bottom}, we can see that errors for $n_X=2976$ in the first column do not improve as $n_Z$ increases. However, in the second column, errors keep reducing as RBF centers are refined. \CMrep{A more fair comparison is to do so diagonally, i.e., we compare}{To make a fair comparison, we can examine the diagonal entries, comparing:}
\begin{itemize}
\item $(n_Z,n_X) = (314,2976)$ and $(744,5296)$ with similar oversampling ratios of 9.5$\times$ and 7.1$\times$, respectively, and
\item $(n_Z,n_X) = (744,2976)$ and $(1320,5296)$ both with oversampling ratios of 4.0$\times$, see \textbf{boldfaced} entries in Table~\ref{tabDC}.
\end{itemize}
In both cases, the reduction in error is trivial.

\begin{figure}
 	\centering
 	\begin{overpic}[width=0.9\textwidth,trim=0 0 0 0]{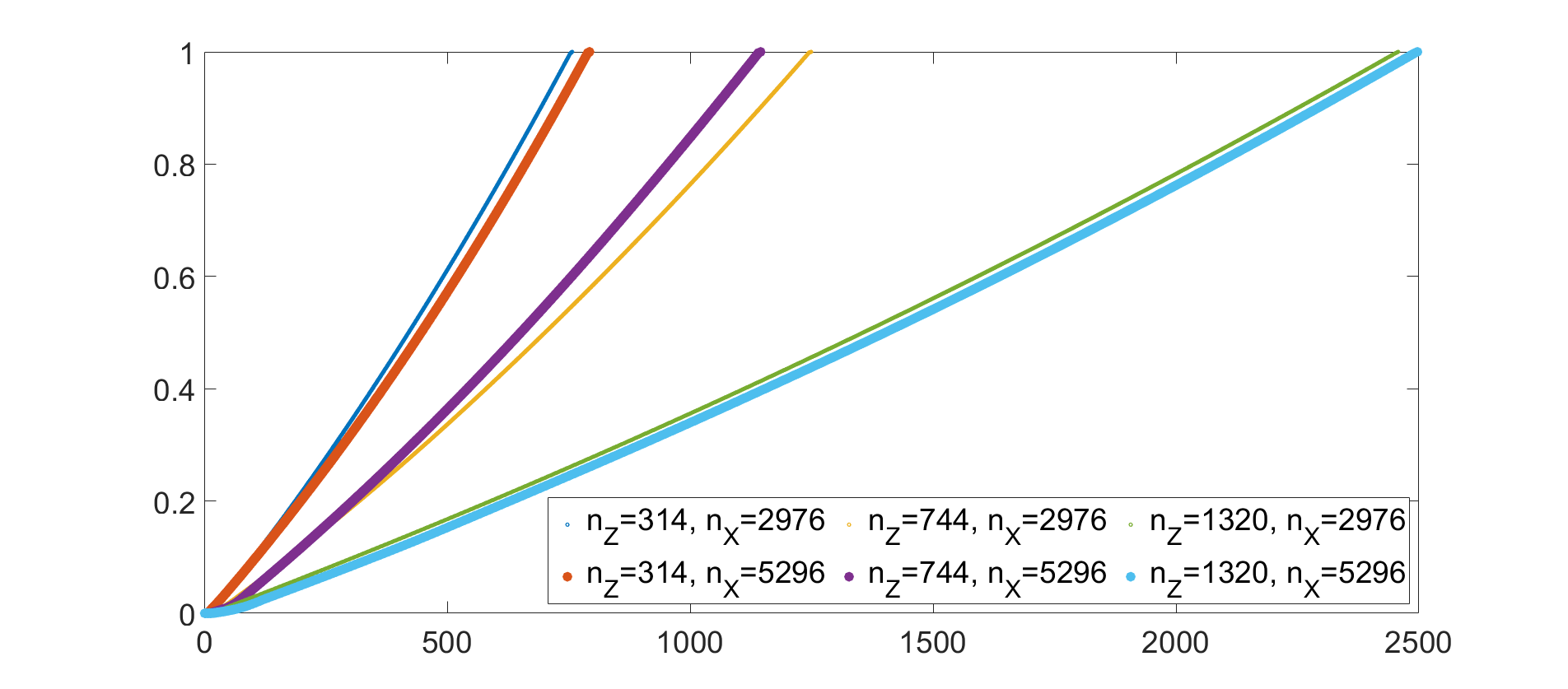}
 		\put(6,18){\scriptsize \rotatebox{90}{{Time $t$}}}
 	\put(50,0){\scriptsize \rotatebox{0}{{\#\,ODE45 Steps}}}
 	\end{overpic}
 	\caption{Exmp.~\ref{ex3}: Current time against the number of adaptive time steps  to complete the integration process  in finding solutions in Figure~\ref{figDC Err}.}\label{figDC dt}
 \end{figure}

We point out that the final accuracy of our results is limited by the ODE45 tolerance. Therefore, all entries in Table~\ref{tabDC} are rather accurate and cannot be used to estimate the convergence rate of the \teal{RBF-LSC-MoL}. Based on these results, we found that using $(n_Z,n_X)=(744,1320)$ provides the best trade-off between accuracy and computation time (measured by ODE45 steps).
To gain some intuition on computational time, we show in Figure~\ref{figDC dt} the number of time steps \teal{required by} ODE45 to solve the problem. In contrast to the fixed $n_Z$ in Figure~\ref{fig4}, we are now seeing results associated with different values of $n_Z$. This figure is easy to summarize: larger $n_Z$ requires \teal{more, hence smaller,} time steps in ODE45, suggesting that the ODE system \eqref{eqQRODE} is getting stiffer with increasing $n_Z$.

 \begin{table}[]
 	\centering
 	\begin{tabular}{|l|llll|}
 	\hline
 		$n_Z \setminus n_X$ & \multicolumn{1}{l|}{744}& \multicolumn{1}{l|}{1320} & \multicolumn{1}{l|}{2976} & \multicolumn{1}{l|}{5296} \\ \hline
 		314        & 1.29e-5 (~~789)& 1.32e-5 (~~789) & 1.19e-5 (~~757) & 8.63e-6 (~~793)              \\ 	\hline
 		744          & NaN          & \blue{8.98e-6 (1333)}  & \textbf{8.91e-6 (1249)}  & 2.80e-6 (1145)             \\ 	\hline
 		1320         & NaN          & NaN             & 8.88e-6 (2457)  & \textbf{2.61e-6 (2497)}             \\ 	\hline
 	\end{tabular}
 	\caption{Exmp.~\ref{ex3}: Maximum errors at $t=1$ (and number of ODE45 time steps taken) of the solutions in Figure~\ref{figDC Err}, see last two columns.
 	}\label{tabDC}
 \end{table}

\bigskip
\exmp{Simulating surface diffusions}{ex4}
 \begin{figure}
 	\centering
 	\begin{overpic}[width=0.5\textwidth]{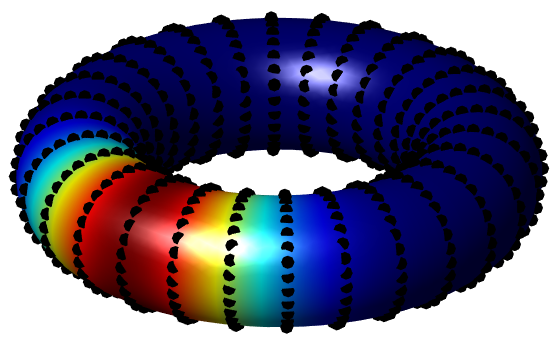}
 	\end{overpic}
 	\caption{Exmp.~\ref{ex4}: Schematic demonstration of data point distributions on Torus with values of the source function $f$ and as colormap.}\label{fig_Torus_2}
 \end{figure}
\begin{figure}
	\centering
	\begin{overpic}[width=0.4\textwidth,trim=0 0 0 0]{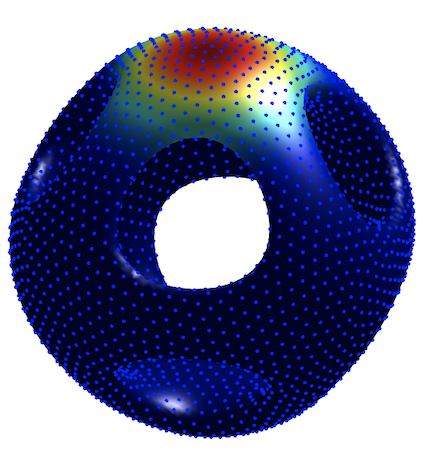}
	\end{overpic}
	\caption{\green{Exmp.~\ref{ex4}: Schematic demonstration of data point distributions on Orthocircle with values of the source function $f$ as colormap}.
}\label{fig_Orthocircle_0}
	\end{figure}
\green{
All observations so far has shown that for achieving the most computationally-efficient setup,
a kernel smoothness of $m=6$ appears to work well with ODE45. Furthermore, it is  beneficial to use a small yet adequate number of RBF centers to ensure enough spatial approximation power. Oversampling by a factor of 1.5 to 2 has been found to be sufficient for obtaining a stable ODE system in our numerical tests.
Exactly how one can pick a proper/optimal kernel and its hyperparameters  is another challenging RBF research topic.
Research into greedy algorithms for RBF interpolation or PDE solving offers some guidance on this issue \cite{Schaback+Wendland-Adapgreetechappr:00,Hon+SchabackETAL-adapgreealgosolv:03,Wenzel-AnalTargDataGree:23,Ling-fastblocalgoquas:16}. These algorithms work by selecting a subset of the provided RBF centers to balance the condition number and accuracy. However, this paper will not go into the specifics of these algorithms.

For this final example, we put our rule-of-thumb into practice}
by simulating a simple physics by considering a surface diffusion
 $\call_\cals=-\Delta_\cals$ in \eqref{equLs} \CMadd{subject to a zero initial condition $g=0$} for $t\in[0,10]$ on a Torus\footnote{Torus: $x^2 + y^2 + z^2 + 1^2-(1/3)^2)^2-4(x^2 + y^2) = 0$} and \CMadd{a Orthocircle\footnote{Orthocircle: $[(x^2 + y^2-1)^2 + z^2][(y^2 + z^2-1)^2 + x^2][(x^2 + z^2-1)^2 + y^2]-0.075^2[1 + 3(x^2 + y^2 + z^2)] = 0.$} respectively. For the Torus, }we can parametrize it with $(x_1,x_2,x_3)\in\M$ as
\[
    (x_1,\,x_2,\,x_3)=\Big(
    \big(1+\cos(\theta)/3\big)\cos(\phi),\,
    \big(1+\cos(\theta)/3\big)\sin(\phi),\,
    \sin(\theta)/3    \Big)
\]
 for $\theta,\phi\in[-\pi,\pi]$.
We use  the right-hand side function $f(\theta,\phi)=\exp({-3\phi^2})$ \CMadd{on the Torus in Figure~\ref{fig_Torus_2} and $f(x,y,z)=\exp(-3((x+0.01)^2+(y-0.003)^2+(z-1.15)^2))$ on the Orthocircle in Figure~\ref{fig_Orthocircle_0}}.

\begin{figure}[]
 	\centering
 	\begin{tabular}{ll}
 		\begin{overpic}[width=0.38\textwidth,trim=55 55 55 55]{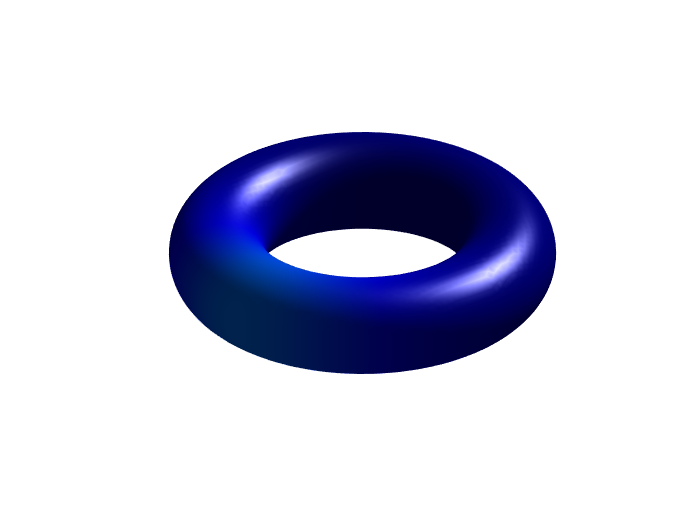}
 			\put(-03,30){\rotatebox{0}{$t=1$}}
 		\end{overpic}
 		\begin{overpic}[width=0.38\textwidth,trim=55 55 55 55]{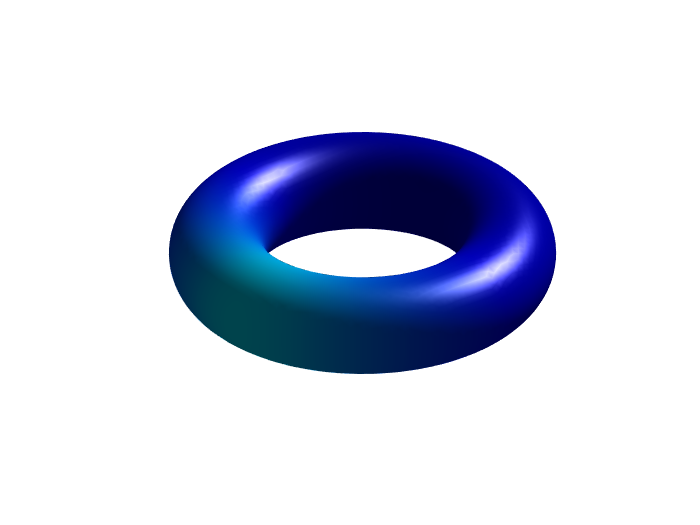}
 			\put(-03,30){\rotatebox{0}{$t=2$}}
 		\end{overpic}
 		\\
 		\begin{overpic}[width=0.38\textwidth,trim=55 55 55 55]{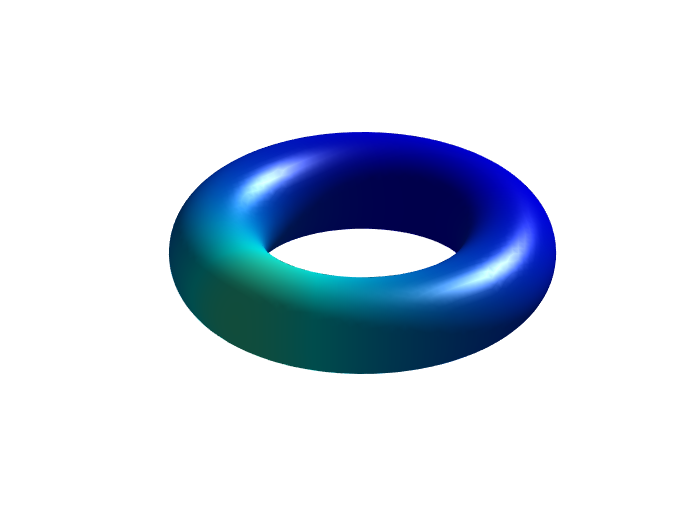}
 			\put(-03,30){\rotatebox{0}{$t=3$}}
 		\end{overpic}
 		\begin{overpic}[width=0.38\textwidth,trim=55 55 55 55]{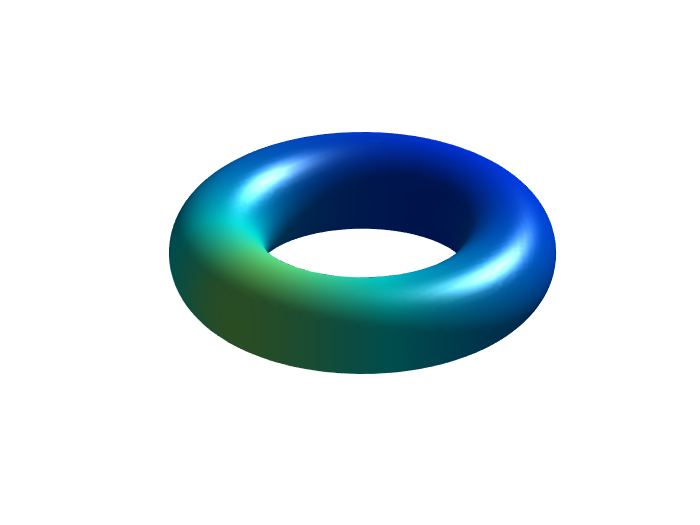}
 			\put(-03,30){\rotatebox{0}{$t=4$}}
 		\end{overpic}
 		\\
 		\begin{overpic}[width=0.38\textwidth,trim=55 55 55 55]{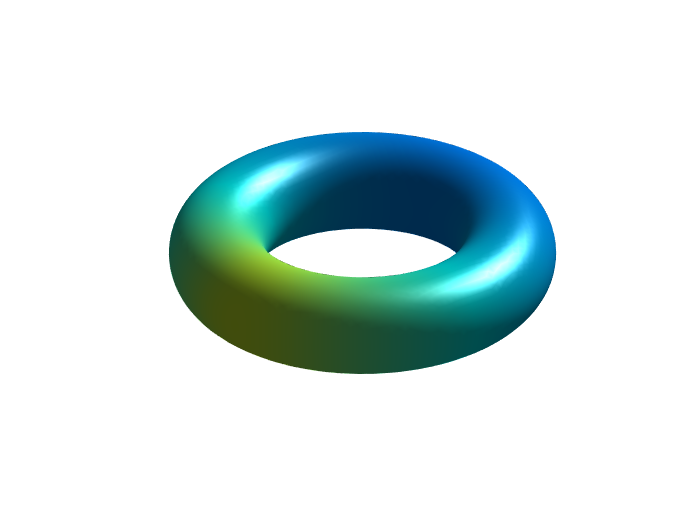}
 			\put(-03,30){\rotatebox{0}{$t=5$}}
 		\end{overpic}
 		\begin{overpic}[width=0.38\textwidth,trim=55 55 55 55]{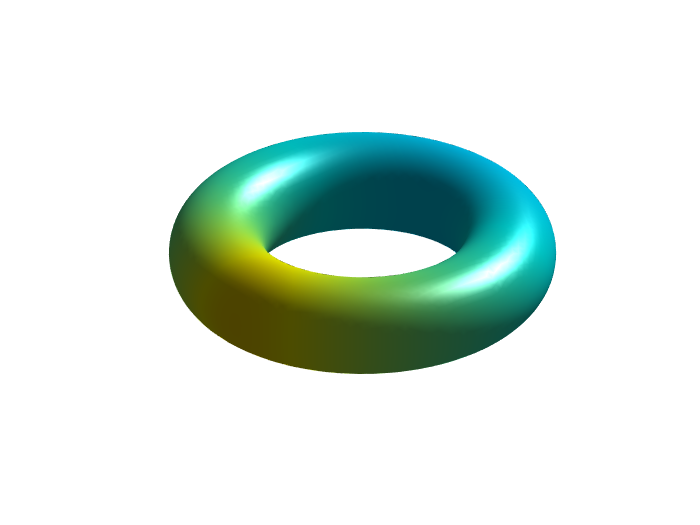}
 			\put(-03,30){\rotatebox{0}{$t=6$}}
 		\end{overpic}
 		\\
 		\begin{overpic}[width=0.38\textwidth,trim=55 55 55 55]{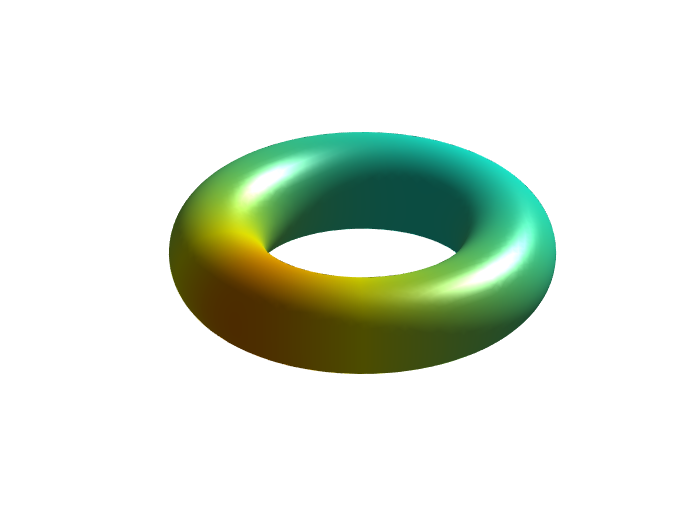}
 			\put(-03,30){\rotatebox{0}{$t=7$}}
 		\end{overpic}
 		\begin{overpic}[width=0.38\textwidth,trim=55 55 55 55]{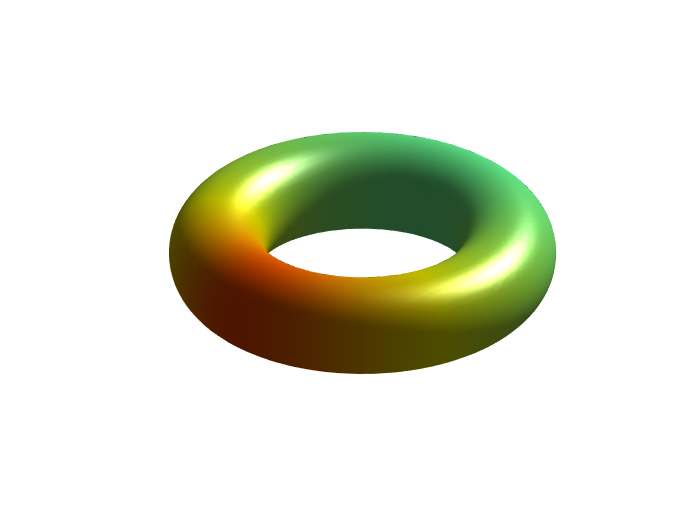}
 			\put(-03,30){\rotatebox{0}{$t=8$}}
 		\end{overpic}
 		\\
 		\begin{overpic}[width=0.38\textwidth,trim=55 55 55 55]{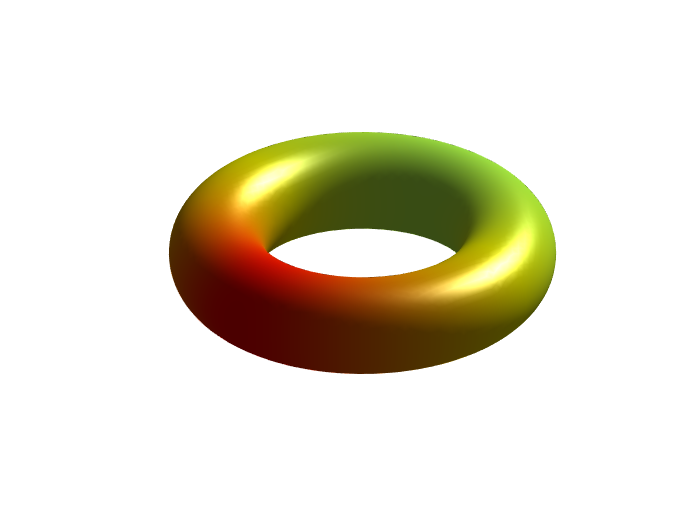}
 			\put(-03,30){\rotatebox{0}{$t=9$}}
 		\end{overpic}
 		\begin{overpic}[width=0.38\textwidth,trim=55 55 55 55]{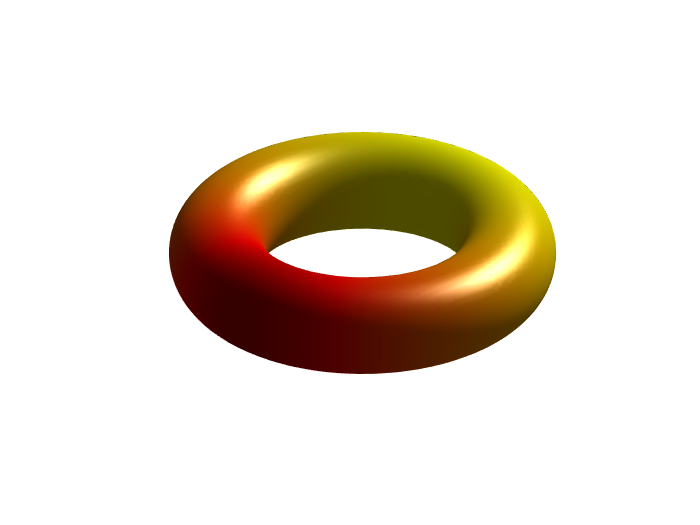}
 			\put(-03,30){\rotatebox{0}{$t=10$}}
 		\end{overpic}
 		\\
 		\multicolumn{2}{r}{\includegraphics[width=0.9\textwidth]{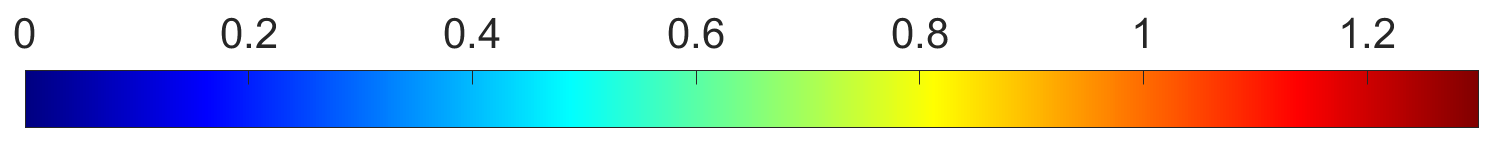}}
 	\end{tabular}
 	\caption{Exmp.~\ref{ex4}: With $n_Z=784$ and $n_X=1156$, snapshots profiles, which use the same colormap for easy comparison, of the RBF method of lines solutions to an isotropic diffusion reaction equation at various time steps for $m=6$ \green{on the Torus}. }\label{figTorus Snap}
 \end{figure}

\begin{figure}[]
	\centering
	\begin{tabular}{ll}
		\begin{overpic}[width=0.35\textwidth,trim=55 30 65 55]{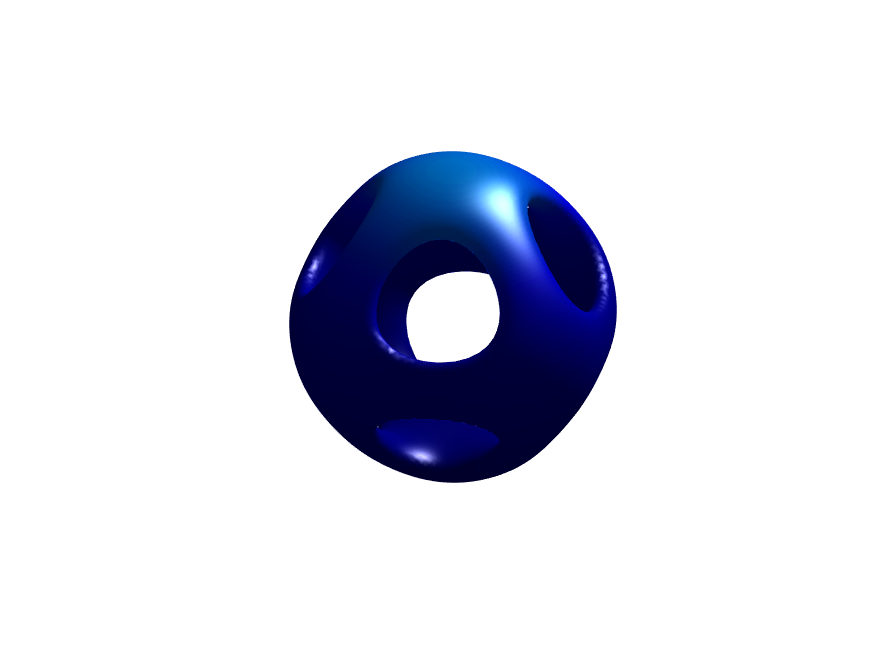}
			\put(-03,40){\rotatebox{0}{$t=1$}}
		\end{overpic}
		\begin{overpic}[width=0.35\textwidth,trim=55 30 65 55]{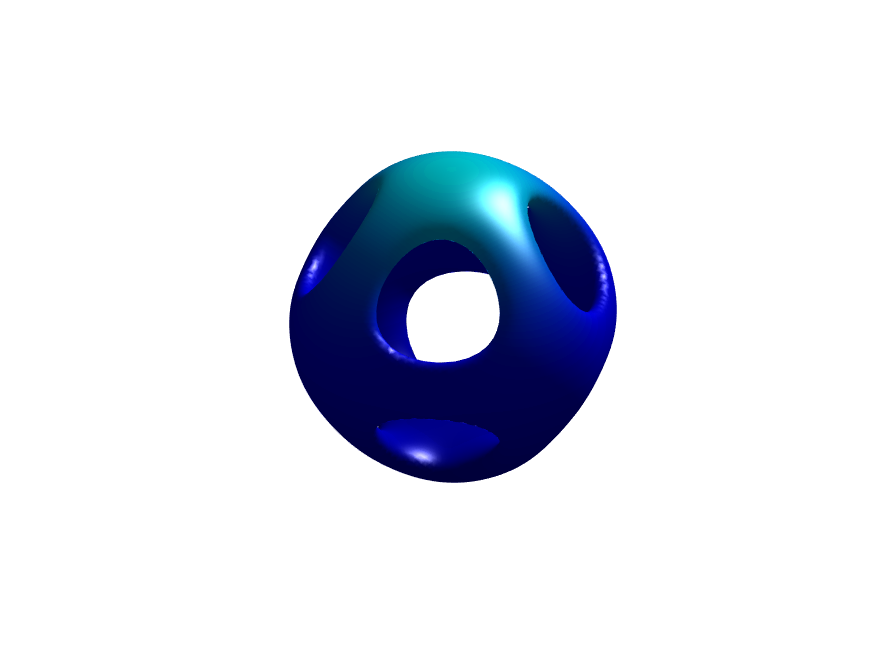}
			\put(-03,40){\rotatebox{0}{$t=2$}}
		\end{overpic}
		\\
		\begin{overpic}[width=0.35\textwidth,trim=55 30 65 55]{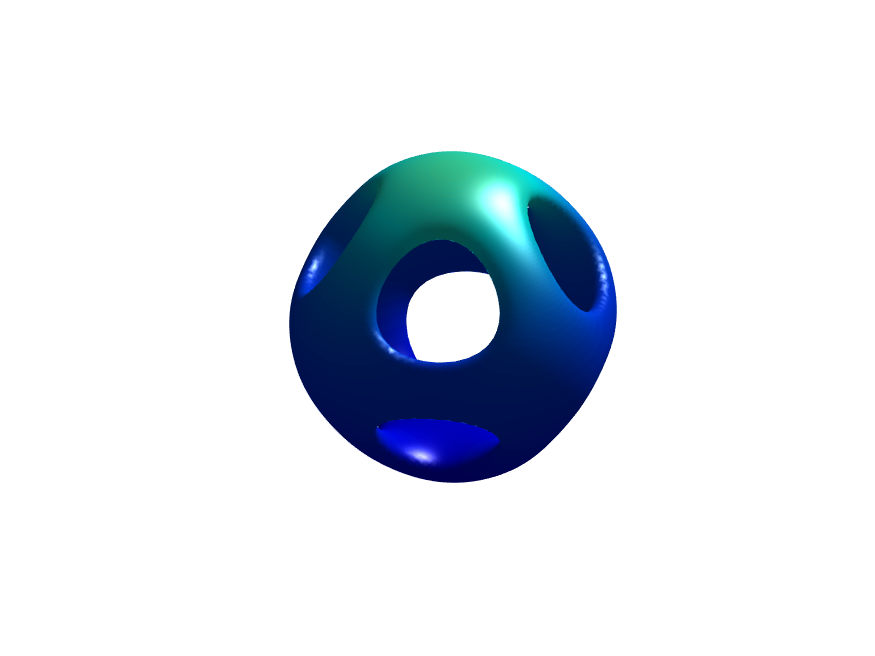}
			\put(-03,40){\rotatebox{0}{$t=3$}}
		\end{overpic}
		\begin{overpic}[width=0.35\textwidth,trim=55 30 65 55]{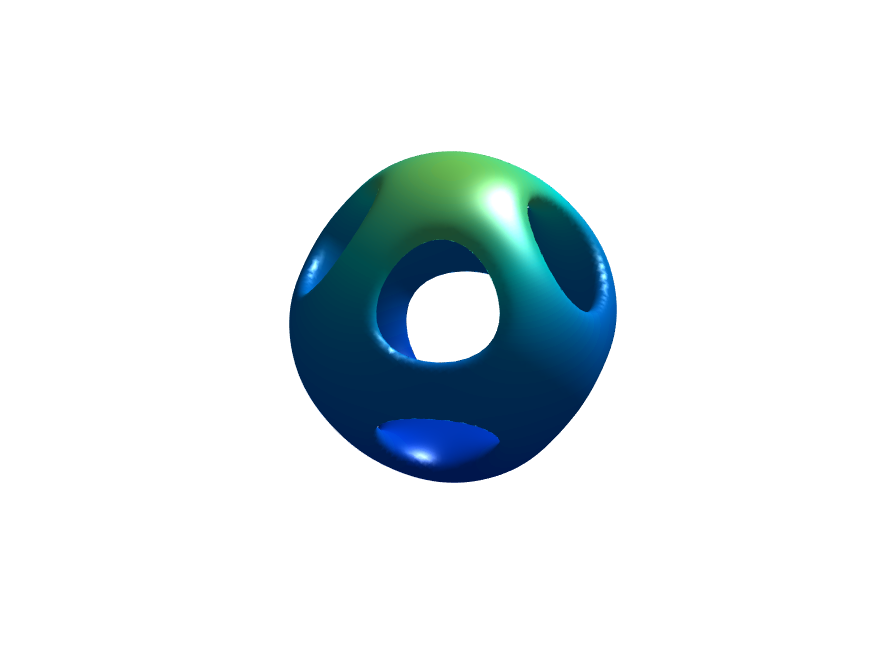}
			\put(-03,40){\rotatebox{0}{$t=4$}}
		\end{overpic}
		\\
		\begin{overpic}[width=0.35\textwidth,trim=55 30 65 55]{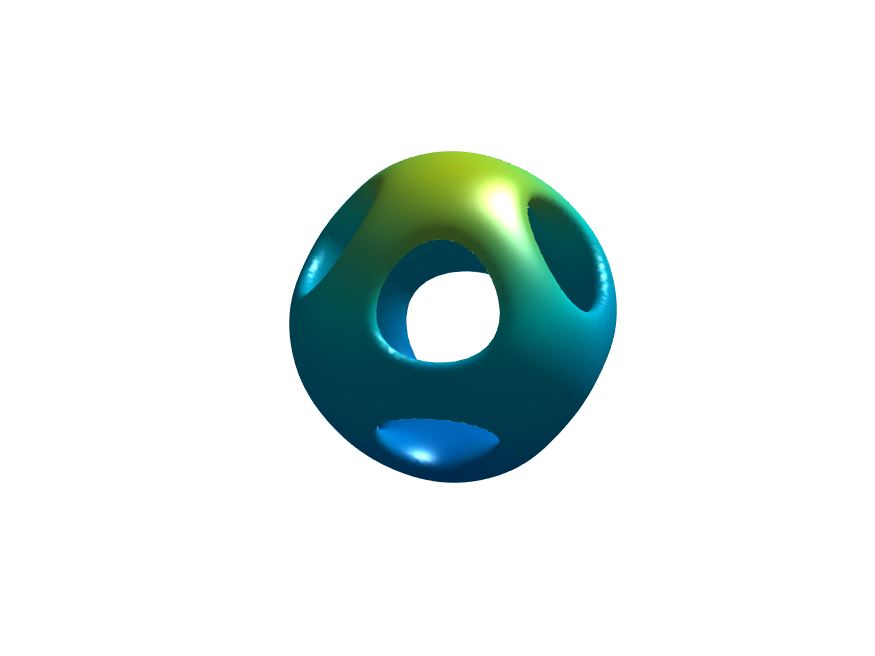}
			\put(-03,40){\rotatebox{0}{$t=5$}}
		\end{overpic}
		\begin{overpic}[width=0.35\textwidth,trim=55 30 65 55]{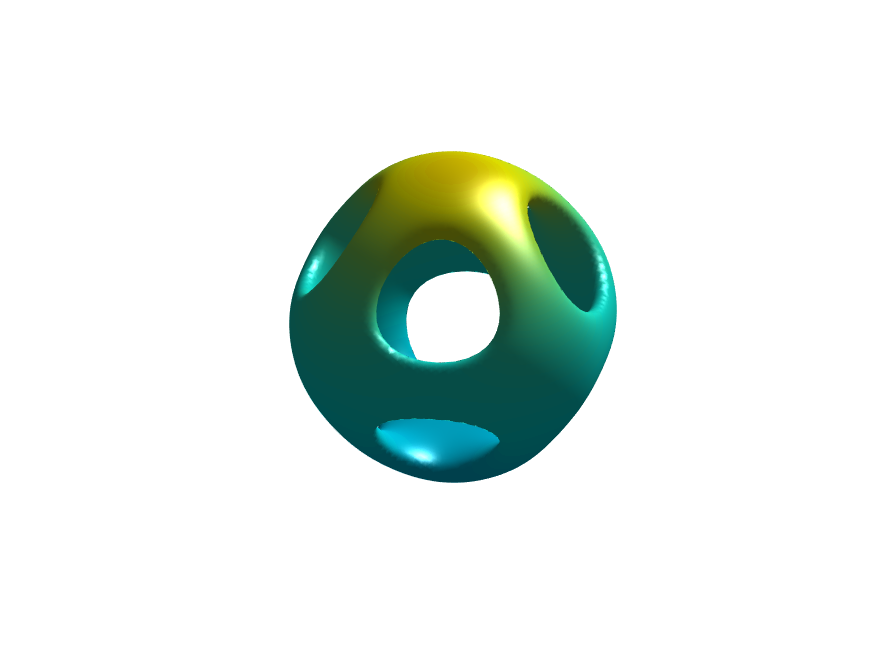}
			\put(-03,40){\rotatebox{0}{$t=6$}}
		\end{overpic}
		\\
		\begin{overpic}[width=0.35\textwidth,trim=55 30 65 55]{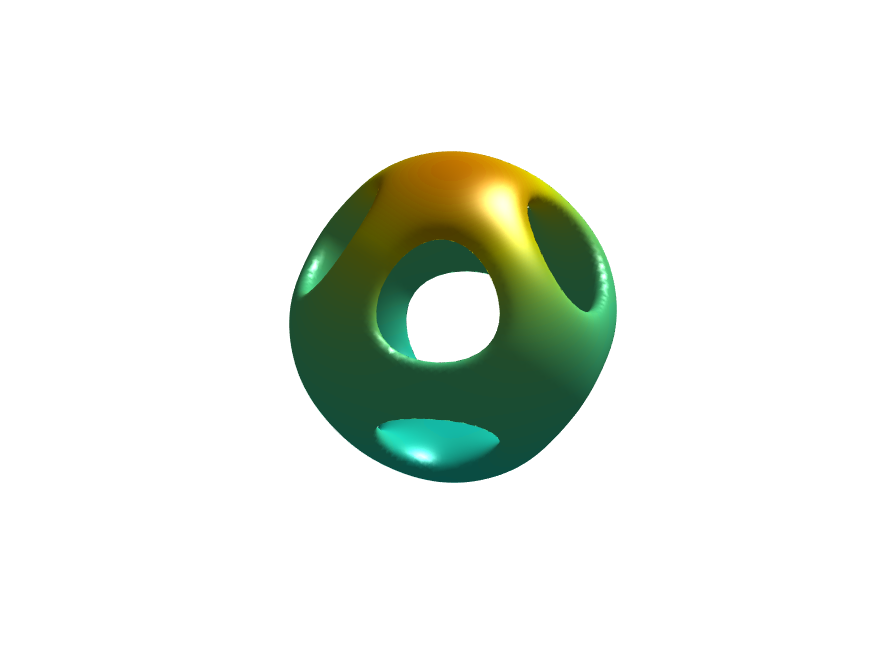}
			\put(-03,40){\rotatebox{0}{$t=7$}}
		\end{overpic}
		\begin{overpic}[width=0.35\textwidth,trim=55 30 65 55]{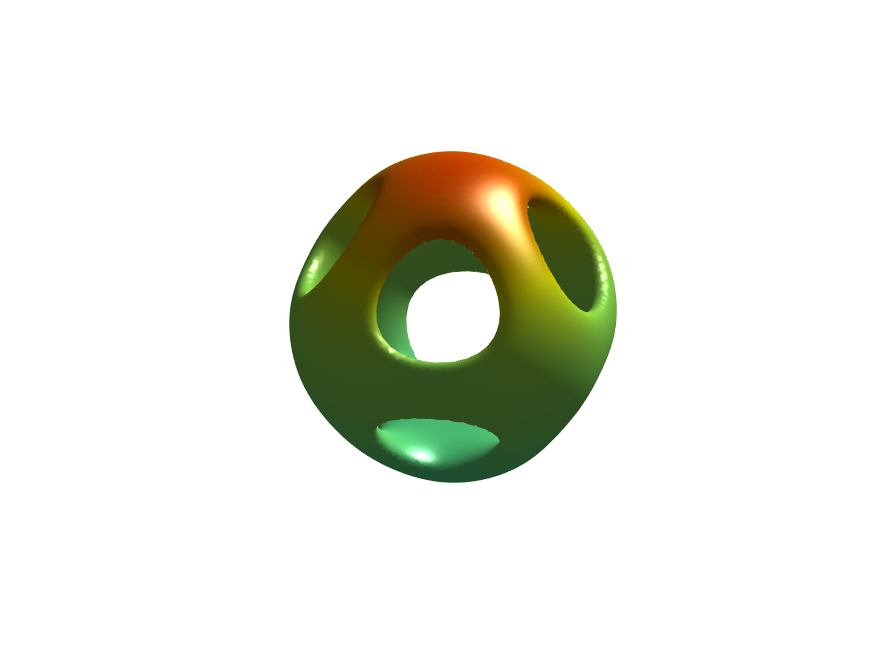}
			\put(-03,40){\rotatebox{0}{$t=8$}}
		\end{overpic}
		\\
		\begin{overpic}[width=0.35\textwidth,trim=55 30 65 55]{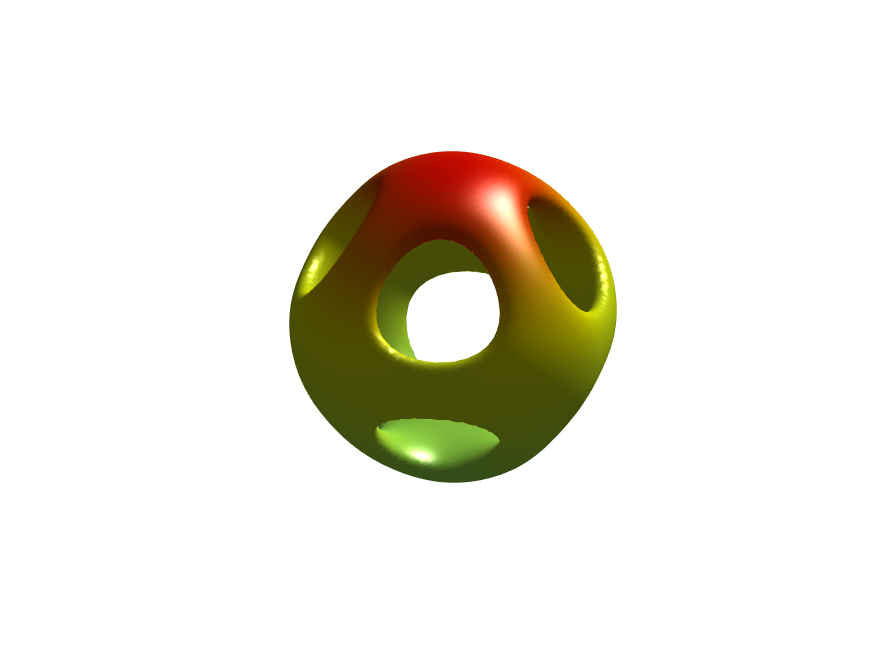}
			\put(-03,40){\rotatebox{0}{$t=9$}}
		\end{overpic}
		\begin{overpic}[width=0.35\textwidth,trim=55 30 65 55]{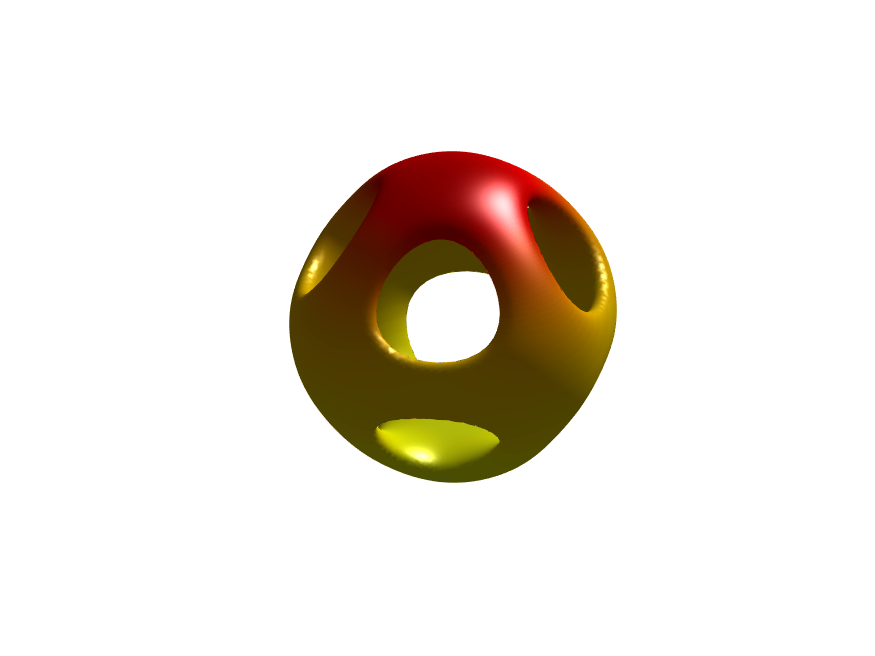}
			\put(-03,40){\rotatebox{0}{$t=10$}}
		\end{overpic}
		\\
		\multicolumn{2}{r}{\begin{overpic}[width=0.9\textwidth,trim=250 1200 150 250]{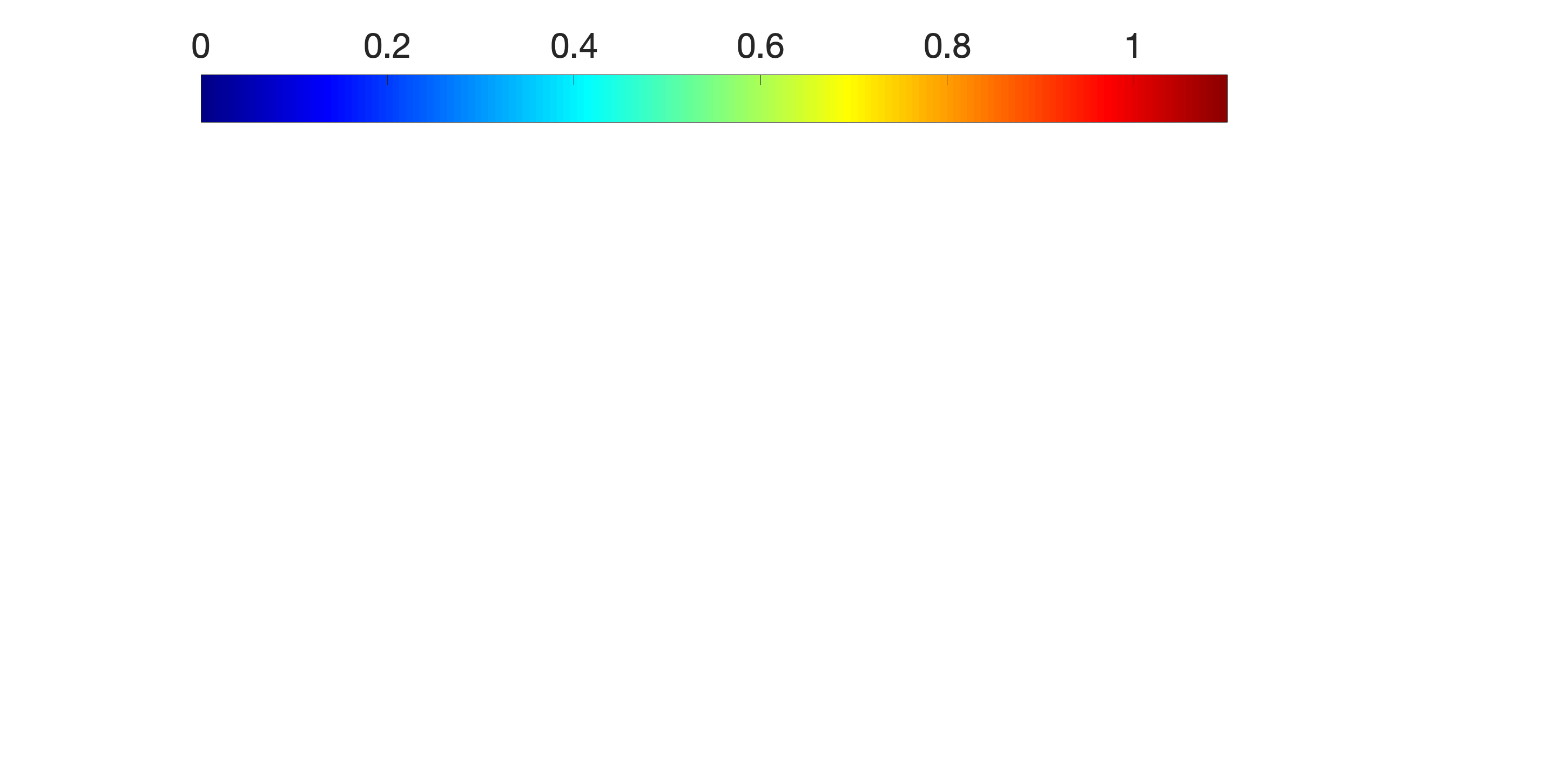}
			\end{overpic}}
	\end{tabular}
	\caption{\green{Exmp.~\ref{ex4}: With $n_Z=3312$ and $n_X=5332$, snapshots profiles, which use the same colormap for easy comparison, of the RBF method of lines solutions to an isotropic diffusion reaction equation at various time steps for $m=4$  on the Orthocircle}. }\label{figOrthocircle Snap}
\end{figure}

\CMadd{On the Torus, we} generate $n_Z=784$ RBF centers and $n_X=1156$ collocation points using the parametric equation given above, as shown in Figure~\ref{fig_Torus_2}. 
The mesh ratios of data points in this example are larger than those in the previous examples.
\CMadd{On the Orthocircle, see  Figure~\ref{fig_Orthocircle_0},  we use quasi-uniform scattered point sets with $n_Z=3312$ and $n_X=5532$, which is the closest to the target 1.5 to 2 oversampling ratio we can get from the quasi-uniform data generator.}
We again apply ODE45 to solve the PDEs and show a few snapshots of the numerical solutions respectively in Figures~\ref{figTorus Snap} \CMadd{and \ref{figOrthocircle Snap}}, from which we can observe the correct physics of heat diffusion from a localized source. Although we do not have analytic solutions to compute errors, we can report that increasing $n_X$ does not result in any noticeable difference in the solutions.

In this paper, we conducted numerical experiments to investigate the effect of oversampling, i.e., using more collocation points than RBF centers, in the RBF least-squares collocation method of lines for solving surface diffusion problems. Our goal was to study the eigenvalue stability in the ODE system after spatial discretization by RBF, rather than targeting high-order convergence. Our results show that oversampling is necessary when using kernels with high smoothness or more RBF basis in the numerical solution. Furthermore, we observed that the number of steps required by the MATLAB function ODE45 increases with the number of RBF centers, but not collocation points.
To achieve the most computationally-efficient setup, we recommend using enough RBF centers that provide sufficient spatial approximation power to reach the temporal error tolerance, and oversampling by \green{a factor of 1.5 to 2.  We successfully run two surface diffusion simulations} based on this rule-of-thumb without fine-tuning.

\bmhead{Acknowledgements}
The authors express their sincere appreciation to the reviewers for their insightful suggestions and thorough review of the manuscript.
 This work was supported by the General Research Fund (GRF No. 12301419, 112301520, 12301021, 12300922) of Hong Kong Research Grant Council, Natural Science Foundation of Jiangxi Province (Grant No. 20212BAB211020), National Natural Science Foundation of China (Grant No. 12001261\green{, 12361086}) and National Key R\&D Program of China (Grant No. 2022YFB4501703).

\bigskip
\bibliographystyle{sn-mathphys} %
\bibliography{Manuscript_NA_review01_LL-minimal}
%
%
%
\end{document}